\documentclass[10pt]{article}
\usepackage[english]{babel}
\usepackage{latexsym}
\usepackage{amsmath}
\usepackage{amssymb}
\usepackage[latin1]{inputenc}
\usepackage{enumerate}
\usepackage{graphicx}
\usepackage{comment}

\usepackage[sc]{mathpazo}

\usepackage{bbm}
\usepackage{mathtools}
\usepackage[amsthm,amsmath,thmmarks,thref]{ntheorem}
\theoremseparator{.}
\newtheorem{definition}{Definition}[section]
\newtheorem{theorem}[definition]{Theorem}
\newtheorem{lemma}[definition]{Lemma}
\newtheorem{remark}[definition]{Remark}
\newtheorem{question}[definition]{Question}
\newtheorem{claim}{Claim}[definition]
\newtheorem{corollary}[definition]{Corollary}
\newtheorem{fact}[definition]{Fact}
\newtheorem{subclaim}[definition]{Subclaim}

\title{The isomorphism relation of theories with S-DOP in generalized Baire spaces}
\author{Miguel Moreno\\ \\University of Helsinki \\Helsinki, Finland\\ \\Bar Ilan University\\Ramat Gan, Israel\\ \\ University of Vienna\\ Vienna, Austria.}
\date{}
\begin{document}
\maketitle
\begin{abstract}
We study the Borel-reducibility of isomorphism relations in the generalized Baire space $\kappa^\kappa$. In the main result we show for inaccessible $\kappa$, that if $T$ is a classifiable theory and $T'$ is superstable with the strong dimensional order property (S-DOP), then the isomorphism of models of $T$ is Borel reducible to the isomorphism of models of $T'$. In fact we show the consistency of the following: If $\kappa$ is inaccessible and $T$ is a superstable theory with S-DOP, then the isomorphism of models of $T$ is $\Sigma_1^1$-complete.
\end{abstract}
\section*{Acknowledgement}
This work was made under the supervision of Tapani Hyttinen. I want to express my gratitude to him for introducing me to the topic, his valuable advices and support during this work. 
I would like to thank the referees for carefully reading the article and providing comments that improve the readability and quality.	
This research was supported by the Doctoral Programme in Mathematics and Statistics of the University of Helsinki (DOMAST). During the revision and editing of this paper, the author was a postdoc at Bar-Ilan University supported by the European Research Council (grant agreement ERC-2018-StG 802756) and was a postdoc at University of Vienna supported by the Vilho, Yrj\"o and Kalle V\"ais\"al\"a Foundation of the Finnish Academy of Science and Letters, and the Austrian Science Fund FWF, Grant I 3709-N35.

\section{Introduction}	

One of the main motivations behind the study of the generalized descriptive set theory, is the connections with model theory. The complexity of a theory can be measured using the Borel reducibility in the generalized Baire spaces: We say that $T'$ is more complex than $T$ if the isomorphism relation of $T$ with universe $\kappa$ ($\cong_T$) is Borel reducible to the isomorphism relation of $T'$ with universe $\kappa$.
Classification theory in Shelah's stability theory gives another notion of complexity. 
The stability theory notion of complexity allows us to compare classifiable theories with non-classifiable theories, but it doesn't allow us to compare the complexity of two non-classifiable theories. On the other hand, the Borel reducibility notion of complexity allows us to compare the complexity of two theories, no matter if the theories are both non-classifiable. Friedman, Hyttinen, Kulikov and others have studied the connection between these two notions of complexity (see \cite{FHK13}).

One of the most important questions regarding the Borel reducibility complexity notion is: \textit{Is the Borel reducibility notion of complexity a refinement of the stability theory notion of complexity?} To answer this question is one of the objectives pursued by the generalized descriptive set theory. For a theory to be non-classifiable, this one must be either unstable, or superstable with OTOP, or superstable with DOP, or stable unsuperstable. It is natural for model theorists to believe that there is a distinction between the complexity of these four kinds of non-classifiable theories, it is conjectured that this may be reflected in the Borel reducibility complexity notion.

The results reviewed in this introduction require further assumptions and the reader is referred to the original paper for the exact assumptions.
In \cite{HKM} it was shown, under the assumptions of $\diamondsuit$ and $\kappa$ successor, if $T$ is classifiable and $T'$ is not, then $\cong_{T}$ is Borel reducible to $\cong_{T'}$. 
In \cite{Fer}, \cite{FMR} and \cite{HKM2} it was showed that for certain models of ZFC, if $\kappa$ is a successor cardinal, then the isomorphism relation of any non-classifiable theory is $\Sigma^1_1$-complete. In particular, in \cite{FMR} and \cite{FMR2} different forcings were constructed to obtain this.
It is natural to ask whether the same holds when $\kappa$ is inaccessible. The case stable unsuperstable was studied in \cite{HM} and the following was found, \textit{if $T$ is classifiable and $T'$ is stable unsuperstable with OCP, then $\cong_T$ is continuously reducible to $\cong_{T'}$}, in the models of \cite{Fer}, \cite{FMR}, and \cite{HKM2} $\cong_{T'}$ is $\Sigma_1^1$-complete.

In \cite{LS} Laskowski and Shelah studied the completeness of $\omega$-stable with DOP. Even though the following question remains open. 
\begin{question}\label{SDOPq}
Is it consistently true: There is a superstable theory with DOP for which the isomorphism relation is $\Sigma^1_1$-complete?
\end{question}

In this article we focus only in cardinals of the form $\kappa^{<\kappa}=\kappa$, the only cases are successors and strong inaccessible.
As it was mentioned above, the previous question was answered when $\kappa$ is a successor (\cite{FMR}, \cite{HKM2}), we will focus only on the case $\kappa$ an inaccessible cardinal. We answer this question in Corollary \ref{C.4.19}, where we show that in the models of \cite{FMR} and \cite{HKM2}, the isomorphism relation of any superstable theory with S-DOP is $\Sigma^1_1$-complete.
In particular we will prove that there is $\lambda<\kappa$ such that $E_{\lambda\text{-club}}^\kappa\leq_c\cong_T$(see Definitions \ref{defequiv} and \ref{D.1.7})  holds for any $T$ superstable theory with S-DOP. In the main result, Corollary \ref{C.4.16}, we settle the Borel reducibility between classifiable theories and superstable theories with S-DOP.

Here and throughout the paper we assume that $\kappa$ is an uncountable
cardinal that satisfies $\kappa^{<\kappa}=\kappa$, $\mathcal{M}$ will denote the monster model, and for every finite tuple $a$, we will denote $a\in A^{length(a)}$ by $a\in A$, unless something else is stated. 

The generalized
Baire space is the set $\kappa^\kappa$ with the bounded topology.  For
every $\zeta\in \kappa^{<\kappa}$, the set 
$$[\zeta]=\{\eta\in \kappa^\kappa \mid \zeta\subset \eta\}$$ 
is a basic open set. The open sets are of the form $\bigcup X$ where
$X$ is a collection of basic open sets. The collection of
Borel subsets of $\kappa^\kappa$ is the smallest set which
contains the basic open sets and is closed under complement and unions of length $\kappa$.

A function $f\colon \kappa^\kappa\rightarrow \kappa^\kappa$ is \emph{Borel}, 
if for every open set $A\subseteq \kappa^\kappa$ the inverse image
$f^{-1}[A]$ is a Borel subset of $\kappa^\kappa$. Let $E_1$ and $E_2$ be
equivalence relations on $\kappa^\kappa$. We say that $E_1$ is 
\emph{Borel reducible} to $E_2$, if there is a Borel function $f$ that satisfies $(x,y)\in E_1\Leftrightarrow
(f(x),f(y))\in E_2$, we call $f$ a \emph{reduction} of $E_1$ to
$E_2$ and it is denoted by $E_1\le_B E_2$. If $f$ is continuous,
then $E_1$ is \emph{continuously reducible} to $E_2$ and it is denoted by $E_1\le_c E_2$.

Let $\mathcal{L}=\{P_n~|~n\in \kappa\backslash\}$ be a given relation vocabulary of size $\kappa$. When we describe a complete theory $T$ in a vocabulary $L\subseteq \mathcal{L}$, we think of it as a complete $\mathcal{L}$-theory extending $T\cup \{\forall \bar{x}\neg
P_n(\bar{x})~|~P_n\in \mathcal{L}\backslash L\}$. We can code $\mathcal{L}$-structures with domain $\kappa$ as follows.

\begin{definition}\label{D.1.6}
  Fix a bijection $\pi\colon \kappa^{<\omega}\to \kappa$. For every
  $\eta\in \kappa^\kappa$ define the $\mathcal{L}$-structure
  $\mathcal{A}_\eta$ with domain $\kappa$ as follows:
  For every relation $P_m$ of arity $n$, every tuple
  $(a_1,a_2,\ldots , a_n)$ in $\kappa^n$ satisfies 
  $$(a_1,a_2,\ldots , a_n)\in P_m^{\mathcal{A}_\eta}\Longleftrightarrow \eta(\pi(m,a_1,a_2,\ldots,a_n))\ge 1.$$
\end{definition}

\begin{definition}[The isomorphism relation]\label{D.1.7} 
Assume $T$ is a complete first order
  theory in a countable vocabulary, $\mathcal{L}$. We define $\cong^\kappa_T$ as the
  relation $$\{(\eta,\xi)\in \kappa^\kappa\times
  \kappa^\kappa\mid (\mathcal{A}_\eta\models T, \mathcal{A}_\xi\models T,
  \mathcal{A}_\eta\cong \mathcal{A}_\xi)\textit{ or }
  (\mathcal{A}_\eta\not\models T, \mathcal{A}_\xi\not\models T)\}.$$
\end{definition}
We will omit the superscript ``$\kappa$'' in $\cong^\kappa_T$ when it
is clear from the context. For every complete first order theory $T$ in a countable
vocabulary there is an isomorphism relation associated with $T$,
$\cong^\kappa_T$.

From now on $\mathcal{L}$ will be a countable relational vocabulary and every theory is a theory in $\mathcal{L}$.

\begin{definition}\label{defequiv}
For every regular cardinal $\mu<\kappa$, $f,g\in \kappa^\kappa$ are $E^\kappa_{\mu\textit{-club}}$ equivalent ($f\  E^\kappa_{\mu\textit{-club}}\ g$) if the set $\{\alpha<\kappa~|~f(\alpha)=g(\alpha)\}$ contains a $\mu$-club, i.e. it is unbounded and closed under $\mu$-limits.
\end{definition}

\section{Preliminaries}\label{sec2}

\subsection{Coloured trees}

Coloured trees have been very useful in the past to reduce $E_{\mu\textit{-club}}^\kappa$ to $\cong_T$ for certain $\mu<\kappa$ and $T$ non-classifiable, see \cite{FHK13}, \cite{HM} or \cite{HK}. We will present a variation of these trees that has height $\lambda+2$ for $\lambda$ an uncountable cardinal.

For a tree $t$, for every $x\in t$ we denote by $ht(x)$ the height of $x$, the order type of $\{y\in t~|~y<x\}$. Define $t_\alpha=\{x\in t~|~ht(x)=\alpha\}$ and $t_{<\alpha}=\cup_{\beta<\alpha}t_\beta$, denote by $x\restriction \alpha$ the unique $y\in t$ such that $y\in t_\alpha$ and $y\leq x$. If $x,y\in t$ and $\{z\in t~|~z<x\}=\{z\in t~|~z<y\}$, then we say that $x$ and $y$ are $\sim$-related, $x\sim y$, and we denote by $[x]$ the equivalence class of $x$ for $\sim$.
An $\alpha, \beta$-tree is a tree $t$ with the following properties:
\begin{itemize}
\item $|[x]|<\alpha$ for every $x\in t$.
\item All the branches have order type less than $\beta$ in $t$.
\item $t$ has a unique root.
\item If $x,y\in t$, $x$ has no immediate predecessors and $x\sim y$, then $x=y$.
\end{itemize}
Notice that in a $\alpha, \beta$-tree, any maximal branch has at most one leaf. Therefore, for any  $\kappa^+$, $(\lambda+2)$-tree we can define functions $c$ over $t_\lambda$. 
\begin{definition}\label{D.2.1}
Let $\lambda$ be an uncountable cardinal. A coloured tree is a pair $(t,c)$, where $t$ is a $\kappa^+$, $(\lambda+2)$-tree and $c$ is a map $c:t_\lambda\rightarrow \kappa\backslash \{0\}$. 
\end{definition}

Two coloured trees $(t,c)$ and $(t',c')$ are isomorphic, if there is a tree isomorphism $f:t\rightarrow t'$ such that for every $x\in t_\lambda$, $c(x)=c'(f(x))$.
We can see every coloured tree as a downward closed subset of $\kappa^{\leq \lambda}$.

For each $f\in \kappa^\kappa$, define the tree $(I_f,d_f)$ as, $I_f$ the set of all strictly increasing functions from some $\theta\leq \lambda$ to $\kappa$ and for each $\eta\in I_f$ with domain $\lambda$, $d_f(\eta)=f(sup(rang(\eta)))$.
For every pair of ordinals $\alpha$ and $\beta$, $\alpha<\beta<\kappa$ and $i<\lambda$ define $$R(\alpha,\beta,i)=\bigcup_{i< j\leq \lambda}\{\eta:[i,j)\rightarrow[\alpha,\beta)~|~\eta \textit{ strictly increasing}\}.$$

Suppose $\kappa$ is an inaccessible cardinal. If $\alpha<\beta<\kappa$ and $\alpha,\beta\neq 0$, let $\{P^{\alpha,\beta}_\gamma~|~\gamma<\kappa\}$ be an enumeration of all downward closed coloured subtrees of $R(\alpha,\beta,i)$ for all $i$, in such a way that each possible coloured tree appears cofinally often in the enumeration. Let $P^{0,0}_0$ be the tree $(I_f,d_f)$.
This enumeration is possible because $\kappa$ is inaccessible; there are at most $|\bigcup_{i<\lambda}\mathcal{P}(R(\alpha,\beta,i))|\leq \lambda\times\kappa=\kappa$ downward closed coloured subtrees, and at most $\kappa\times \kappa^{<\kappa}=\kappa$ coloured trees.
Denote by $Q(P^{\alpha,\beta}_\gamma)$ the unique ordinal number $i$ such that $P^{\alpha,\beta}_\gamma\subset R(\alpha,\beta,i)$.
\begin{definition}\label{D.2.6}
Suppose $\kappa$ is an inaccessible cardinal. 
Order the set $\lambda\times \kappa\times \kappa\times \kappa\times \kappa$ lexicographically, $(\alpha_1,\alpha_2,\alpha_3,\alpha_4,\alpha_5)>(\beta_1,\beta_2,\beta_3,\beta_4,\beta_5)$ if for some $1\leq k \leq 5$, $\alpha_k>\beta_k$ and for every $i<k$, $\alpha_i=\beta_i$. Order the set $(\lambda\times \kappa\times \kappa\times \kappa\times \kappa)^{\leq \lambda}$ as a tree by inclusion.
For every $f\in \kappa^\kappa$ define $J_f=(J_f,c_f)$ as the subtree of all $\eta: s\rightarrow \lambda\times \kappa^4$, where $s\leq \lambda$, ordered by extension, and such that the following conditions hold for all $i,j<s$:

Denote by $\eta_i$, $1\leq i\leq 5$, the functions from $s$ to $\kappa$ that satisfy, $\eta(n)=(\eta_1(n),\eta_2(n),\eta_3(n),\eta_4(n),\eta_5(n))$.
We will construct the coloured tree $J_f$ intuitively by merging many different coloured trees, $\eta_1$ fixes a domain of the form $[i,j)$, the functions $\eta_2$,  $\eta_3$, and $\eta_4$ fix a coloured tree of the form $P^{\alpha,\beta}_{\gamma}\subset R(\alpha,\beta,i)$, and $\eta_5$ is an element of $P^{\alpha,\beta}_{\gamma}$ that is responsible to define the color. 
\begin{enumerate}
\item [1.] $\eta\restriction n\in J_f$ for all $n<s$.
\item [2.]$\eta$ is strictly increasing with respect to the lexicographical order on $\lambda\times \kappa^4$.
\item [3.]$\eta_1(i)\leq \eta_1(i+1)\leq \eta_1(i)+1$.
\item [4.]$\eta_1(i)=0$ implies $\eta_2(i)=\eta_3(i)=\eta_4(i)=0$.
\item [5.]$\eta_2(i)\ge\eta_3(i)$ implies $\eta_2(i)=0$.
\item [6.]$\eta_1(i)<\eta_1(i+1)$ implies $\eta_2(i+1)\ge \eta_3(i)+\eta_4(i)$.
\item [7.]For every limit ordinal $\alpha$, $\eta_k(\alpha)=sup_{\beta<\alpha}\{\eta_k(\beta)\}$ for $k\in \{1,2\}$.
\item [8.]$\eta_1(i)=\eta_1 (j)$ implies $\eta_k (i)=\eta_k (j)$ for $k\in \{2,3,4\}$.
\item [9.]If for some $k<\lambda$, $[i,j)=\eta_1^{-1}\{k\}$, then $$\eta_5\restriction {[i,j)}\in P^{\eta_2(i),\eta_3(i)}_{\eta_4(i)}.$$
\noindent Note that 7 implies $Q(P^{\eta_2(i),\eta_3(i)}_{\eta_4(i)})=i$.
\item [10.]If $s=\lambda$, then either 
\begin{itemize}
\item [(a)] there exists an ordinal number $m$ such that for every $k<m$ $\eta_1(k)<\eta_1(m)$, for every $k \ge m$ $\eta_1(k)=\eta_1(m)$, and the color of $\eta$ is determined by $P^{\eta_2(m),\eta_3(m)}_{\eta_4(m)}$: $$c_f(\eta)=c(\eta_5\restriction {[m,\lambda)})$$ where $c$ is the colouring function of $P^{\eta_2(m),\eta_3(m)}_{\eta_4(m)}$.
\end{itemize}
Or
\begin{itemize}
\item [(b)] there is no such ordinal $m$ and then $c_f(\eta)=f(sup(rang(\eta_5)))$.
\end{itemize}
\end{enumerate}
\end{definition}
The following lemma is a variation of Lemma 4.7 of \cite{HM}, 
nevertheless the proof is the same in both cases.
\begin{lemma}\label{L.2.7}
Assume $\kappa$ is an inaccessible cardinal, then for every $f,g\in \kappa^\kappa$ the following holds $$f\  E^\kappa_{\lambda\textit{-club}}\ g \Leftrightarrow  J_f\cong J_g$$
\end{lemma}
\begin{remark}\label{R.2.7a}
For each $\alpha<\kappa$ define $J_f^\alpha$ as $$J_f^\alpha=\{\eta\in J_f~|~ rang(\eta)\subset \lambda\times(\beta)^4\textit{ for some }\beta<\alpha\}.$$ Notice that every $\eta\in J_f$ has the following properties:
\begin{enumerate}
\item [1.]
$
sup(rang(\eta_4))\leq sup(rang(\eta_3))=sup(rang(\eta_5))=sup(rang(\eta_2)).
$
\item [2.]When $\eta\restriction k\in J_f^\alpha$ holds for every $k\in \lambda$, $sup(rang(\eta_5))\leq \alpha$. If in addition $\eta\notin J_f^\alpha$, then 
$
sup(rang(\eta_5))= \alpha.
$
\end{enumerate} 
\end{remark}
Let us take a look at the sets $rang(f)$ and $rang(c_f)$, more specifically at the set $\{\alpha<\kappa~|~f(\alpha)\in rang(c_f)\}$.
\begin{remark}\label{R.2.8}
Assume $f\in \kappa^\kappa$ and let $J_f$ be the respective coloured tree obtained by Definition \ref{D.2.6}. If $\eta\in J_f$ satisfies Definition \ref{D.2.6} item 10 (b), then clearly there exists $\alpha<\kappa$ such that $c_f(\eta)=f(\alpha)$. It is possible that not for every $\alpha<\kappa$, there is $\eta\in J^{\alpha+1}_f$ such that $c_f(\eta)=f(\alpha)$. Nevertheless the set $C=\{\alpha<\kappa~|~\exists \xi\in J^{\alpha+1}_f\textit{ such that }\xi_1\restriction \omega=id+1, \xi_1\restriction {[\omega,\lambda)}=id\restriction {[\omega,\lambda)}\textit{ and } c_f(\xi)=f(\alpha)\}$ is a $\lambda$-club. 
\end{remark}

\subsection{Strong DOP}

Now, we will recall the dimensional order property and the strong dimensional order property. The independence properties of indiscernible sequences have been a very useful tool to study theories with DOP (see \cite{HaMa}, Section 2), this makes superstable theories with DOP and strong independence properties good candidates to answer Question \ref{SDOPq}. Following this direction we will define the strong dimensional property (Lemma \ref{L.3.9} and Definition \ref{D.3.13}), we will give some important properties that will be useful to construct models of theories with the strong dimensional property. It is an open question whether every superstable theory with DOP has also S-DOP

In \cite{Sh} Shelah gives an axiomatic approach for an isolation notion, $F$, and defines the notions $F$-constructible, $F$-atomic, $F$-primary, $F$-prime and $F$-saturated. 

\begin{definition}\label{D.3.1}
Denote by $F_\theta^a$ the set of pairs $(p,B)$ with $|B|<\theta$, such that for some $A\supseteq B$ and $c$, $p\in S(A)$, $c\models p$ and $stp(c,B)\vdash p$.
\end{definition}
In \cite{Sh} (Definition III 4.2 (2), and Definition V 1.1 (2) and (4)) the notions of stationarization of a type, and orthogonal types are defined.
For $p_1,p_2\in S(A)$ stationary types the following holds. If $p_1=tp(a_1,A)$, and $p_2=tp(a_2,A)$, then $p_1$ is weakly orthogonal to $p_2$ if and only if $a_1\downarrow_A a_2$. 
A stationary type $p\in S(B)$ is orthogonal to $A$ if for all $a,b$ and $D\supset A$ the following holds: If $tp(b,B)$ is stationary, $a\models p$, $b\downarrow_AB$, $b\downarrow_BD$ and $a\downarrow_BD$, then $a\downarrow_Db$.
\begin{fact}\label{F.3.5}
Let $B,D\subseteq M$, $M$ a $F_\omega^a$-saturated model over $B\cup D$, and $p\in S(M)$. If $p$ is orthogonal to $D$ and $p$ does not fork over $B\cup D$, then for every $c\models p\restriction B\cup D$ the following holds: $c\downarrow_{B\cup D}M$ implies $tp(c,M)\perp D$.
\end{fact}

A type $p\in S(B\cup C)$ is orthogonal to $C$, if for every $F_\omega^a$-primary model, $M$, over $B\cup C$ there exists a non-forking extension of $p$, $q\in S(M)$, orthogonal to $C$.

In \cite{Sh} (X.2 Definition 2.1) Shelah defines the dimensional order property, DOP, as follows.
\begin{definition}\label{D.3.7}
A theory $T$ has the dimensional order property (DOP) if there are $F_{\kappa(T)}^a$-saturated models $(M_i)_{i<3}$, $M_0\subset M_1\cap M_2$, $M_1\downarrow_{M_0}M_2$, and the $F_{\kappa(T)}^a$-prime model, $M_3$, over $M_1\cup M_2$ is not $F_{\kappa(T)}^a$-minimal over $M_1\cup M_2$.
\end{definition} 

\begin{definition}
Let $I$ be an infinite indiscernible set. We define $Av(I,A)$ the average type of $I$ over $A$, to be the set $$\{\varphi(x,a)\mid a\in A, |\{b\in I\mid\ \mathcal{M} \models  \varphi(b,a)\}|\ge \omega\}.$$
\end{definition}

Recall that an $F_{\kappa(T)}^a$-prime model is $F_{\kappa(T)}^a$-atomic and $F_{\kappa(T)}^a$-saturated.
The rest of the results in this section will be stated and proved for the case of the $F_{\omega}^a$ isolation. Many of those results can be easily generalized to $F_{\kappa(T)}^a$ by making small changes on the proof.
From now on we will work only with superstable theories. We know that for every superstable theory $T$, $\kappa(T)=\omega$.
The proof of the following lemma is similar to the proof of [\cite{Sh} X.2 Lemma 2.2].

\begin{lemma}\label{L.3.9}
Let $M_0\subset M_1\cap M_2$ be $F_{\omega}^a$-saturated models, $M_1\downarrow_{M_0}M_2$, $M_3$ $F_{\omega}^a$-atomic over $M_1\cup M_2$ and $F_{\omega}^a$-saturated. Then the following conditions are equivalent:
\begin{enumerate}
\item [1.]There is a non-algebraic type $p\in S(M_3)$ orthogonal to $M_1$ and to $M_2$, that does not fork over $M_1\cup M_2$.
\item [2.]There is an infinite indiscernible $I\subseteq M_3$ over $M_1\cup M_2$ that is independent over $M_1\cup M_2$.
\item [3.]There is an infinite $I\subseteq M_3$ indiscernible over $M_1\cup M_2$ and independent over $M_1\cup M_2$, such that $Av(I,M_3)$ is orthogonal to $M_1$ and to $M_2$.
\end{enumerate}
\end{lemma}

\begin{lemma}[\cite{HS99}, Theorem 2.1]\label{L.3.10}
Let $M_0\prec M_1,M_2$ be $F_{\omega}^a$-saturated models, such that $M_1\downarrow_{M_0}M_2$. Let $M_3$ be an $F_{\omega}^a$-prime model over $M_1\cup M_2$ and let $I\subseteq M_3$ be an indiscernible over $M_1\cup M_2$ such that $Av(I,M_3)$ is orthogonal to $M_1$ and to $M_2$. If $(B_i)_{i<3}$ are sets such that: 
\begin{itemize}
\item $B_0\downarrow_{M_0}M_1\cup M_2$.
\item $B_1\downarrow_{M_1\cup B_0}B_2\cup M_2$.
\item $B_2\downarrow_{M_2\cup B_0}B_1\cup M_1$.
\end{itemize}
Then $$tp(I,M_1\cup M_2) \vdash tp(I,M_1\cup M_2\cup_{i<3}B_i).$$
\end{lemma}

The following lemma give us sufficient conditions over models $M_1'$, $M_2'$, and $M_3'$ such that, if $M_1$, $M_2$, and $M_3$ are models that satisfy Definition \ref{D.3.7}, then $M_1'$, $M_2'$, and $M_3'$ satisfy Definition \ref{D.3.7}.

\begin{lemma}\label{L.3.11}
Let $M_0\subset M_1\cap M_2$ be $F_{\omega}^a$-saturated models, such that $M_1\downarrow_{M_0}M_2$ and $M_3$, the $F_{\omega}^a$-prime model over $M_1\cup M_2$, is not $F_{\omega}^a$-minimal over $M_1\cup M_2$.
If $(M'_i)_{i<3}$ are $F_{\omega}^a$-saturated models that satisfy:
\begin{itemize}
\item $\forall i<3$, $M_i\subseteq M'_i$.
\item $\forall i<3$, $M'_i\downarrow_{M_i}M_3$.
\item $M'_1\downarrow_{M'_0}M'_2$.
\end{itemize}
Then $M_3'$ the $F_{\omega}^a$-prime model over $M'_1\cup M_2'$ is not $F_{\omega}^a$-minimal over $M_1'\cup M_2'$.
\end{lemma}

\begin{remark}\label{R.3.12}
From the previous lemma we can conclude that if $I$ is independent over $M_1\cup M_2$, then 
$I$ is independent over $M'_1\cup M'_2$.
\end{remark}
\begin{definition}\label{D.3.13}
We say that a superstable theory $T$ has the strong dimensional order property (S-DOP) if the following holds:

There are $F_{\omega}^a$-saturated models $(M_i)_{i<3}$, $M_0\subset M_1\cap M_2$, such that $M_1\downarrow_{M_0}M_2$, and for every $M_3$ $F_{\omega}^a$-prime model over $M_1\cup M_2$, there is a non-algebraic type $p\in S(M_3)$ orthogonal to $M_1$ and to $M_2$, such that it does not fork over $M_1\cup M_2$.
\end{definition}

From \cite{Sh} X.2 Lemma 2.2, every superstable theory with S-DOP has DOP.
In \cite{HrSo} Hrushovski and Sokolvi\'c proved that the theory of differentially closed fields of characteristic zero (DCF) has eni-DOP, so it has DOP. The reader can find an outline of this proof in \cite{Mar07}. We will show that DCF also has the S-DOP, this can be done following the proof in \cite{Mar07} or the one in \cite{Mar} which uses Rosenlicht's Theorem. 

More on DCF (proofs, definitions, references, etc) can be found in \cite{Mar}.
Let $K$ be a saturated model of DCF, $k\subseteq K$ and $a\in K^n$, we denote by $k\langle a\rangle$ the differentially closed subfield generated by $k(a)$. If $A\subseteq K$ and for all $n$, every nonzero $f\in k\{x_1,x_2,\ldots,x_n\}$, and all $a_1,a_2,\ldots,a_n\in A$ it holds that $f(a_1,a_2,\ldots,a_n)\neq 0$, then we say that $A$ is $\delta$-independent over $k$. 
For all $k\subseteq K$ denote by $k^{dif}$ the differential closure of $k$ in $K$.

\begin{theorem}[Hrushovski, Sokolvi\'c, \cite{Mar} Theorem 7.6, \cite{Mar07} Sections 4, 5 ]\label{T.3.17}
Suppose $K_0$ is a differentially closed field with characteristic zero, $\{a,b\}$ is $\delta$-independent over $K_0$, $K_1=K_0\langle a\rangle^{dif}$, $K_2=K_0\langle b\rangle^{dif}$, and $K=K_0\langle a,b\rangle^{dif}$.
 There is a type $p$ over $K$ that does not fork over $\{a,b\}$ such that $K_1\downarrow_{K_0}K_2$, $p\perp K_1$, and $p\perp K_2$.
\end{theorem}
\begin{corollary}\label{C.3.18}
DCF has the S-DOP.
\end{corollary}
\begin{proof}
Let $a$, $b$, $K_1$, $K_2$, and $p$ be as in Theorem \ref{T.3.17}. By Theorem \ref{T.3.17} it is enough to show that $p$ does not fork over $K_1\cup K_2$. This follows since $p$ does not fork over $\{a,b\}$. 
\end{proof}

\section{Construction of models}

\underline{From now on $\kappa$ will be an inaccessible cardinal.}

In this section we will use coloured trees to construct models of a superstable theory with S-DOP. To do this, we will need some basic results first and fix some notation. Instead of writing $F_\omega^a$-constructible, $F_\omega^a$-atomic, $F_\omega^a$-saturated and $F_\omega^a$-saturated we will write $a$-constructible, $a$-atomic, $a$-primary, $a$-prime and $a$-saturated. From now on $T$ will be a superstable theory with S-DOP.
We will denote by $\lambda(T)$ the least cardinal such that $T$ is $\lambda$-stable. Since $T$ is superstable, then $\lambda(T)\leq 2^\omega$, we will denote by $\lambda$ the cardinal $(2^\omega)^+$.
\begin{definition}\label{D.4.1}
\begin{itemize}
\item Let us define the dimension of a type $p\in S(A)$ in $M$ by: $dim(p,M)=min\{|J|:J\subseteq M$, $J$ is a maximal independent sequence over $A$, and $\forall a \in J, a\models p\}$
\item Let us define the dimension of an indiscernible $I$ over $A$ in $M$ by:
$dim(I,A,M)=min\{|J|:J$ is equivalent to $I$ and $J$ is a maximal indiscernible over $A$ in $M\}$. If for all $J$ as above $dim(I,A,M)=|J|$, then we say that the dimension is true.
\end{itemize}
\end{definition}
\begin{lemma}[\cite{Sh}]\label{L.4.2}
If $I$ is a maximal indiscernible set over $A$ in $M$, then $|I|+\kappa(T)=dim(I,A,M)+\kappa(T)$, 
and if $dim(I,A,M)\ge\kappa(T)$, then the dimension is true.
\end{lemma}
\begin{theorem}[\cite{Sh}]\label{T.4.3}
If $M$ is an $a$-primary model over $A$, and $I\subseteq M$ is an infinite indiscernible set over $A$, then $dim(I,A,M)=\omega$.
\end{theorem}
For any indiscernible sequence $I=\{a_i~|~i<\gamma\}$, we will denote by $I\restriction _\alpha$ the sequence $I=\{a_i~|~i<\alpha\}$. Now for every $f\in \kappa^\kappa$ we will use the the tree $J_f$ given in Definition \ref{D.2.6}, to construct the model $\mathcal{A}^f$.
Since $T$ has the S-DOP, by Lemma \ref{L.3.9} and Lemma \ref{L.3.10} there are $a$-saturated models $\mathcal{A}, \mathcal{B}, \mathcal{C}$ of cardinality $2^\omega$
and an indiscernible sequence $\mathcal{I}$ over $\mathcal{B}\cup \mathcal{C}$ of size $\kappa$ that is independent over $\mathcal{B}\cup \mathcal{C}$ such that 
\begin{enumerate}
\item [1.]$\mathcal{A}\subset\mathcal{B}\cap\mathcal{C}$, $\mathcal{B}\downarrow_{\mathcal{A}}\mathcal{C}$.
\item [2.]$Av(\mathcal{I},\mathcal{B}\cup\mathcal{C})$ is orthogonal to $\mathcal{B}$ and to $\mathcal{C}$.
\item [3.]If $(B_i)_{i<3}$ are sets such that:
\begin{itemize}
\item [(a)]$B_0\downarrow_{\mathcal{A}}\mathcal{B}\cup\mathcal{C}$.
\item [(b)]$B_1\downarrow_{\mathcal{B}\cup B_0}B_2\cup\mathcal{C}$.
\item [(c)]$B_2\downarrow_{\mathcal{C}\cup B_0}B_1\cup\mathcal{B}$.
\end{itemize}
Then, $$tp(\mathcal{I},\mathcal{B}\cup\mathcal{C}) \vdash tp(\mathcal{I},\mathcal{B}\cup\mathcal{C}\cup_{i<3}B_i).$$
\end{enumerate}
By the existence property of forking, for any $D\supseteq \mathcal{A}$ there is $F\in Aut(\mathcal{A})$ such that  for all $c\in \mathcal{B}$, $stp(F(c),\mathcal{A})=stp(c, \mathcal{A})$ and $F(c)\downarrow_\mathcal{A} D$ (the same holds for $\mathcal{C}$).
For every $\xi\in (J_f)_{<\lambda}$ and every $\eta\in (J_f)_{\lambda}$ (recall $t_\alpha$ at the beginning of the section \ref{sec2}), let $\mathcal{B}_\xi\cong_\mathcal{A}\mathcal{B}$, $\mathcal{A}\preceq \mathcal{B}_\xi$, and
 $\mathcal{C}_\eta\cong_\mathcal{A}\mathcal{C}$, $\mathcal{A}\preceq \mathcal{C}_\eta$, such that the models $(\mathcal{B}_\xi)_{\xi\in (J_f)_{<\lambda}}$ and $(\mathcal{C}_\eta)_{\eta\in (J_f)_{\lambda}}$ satisfy the following:
 \begin{itemize}
 \item $\mathcal{B}_\xi\downarrow_\mathcal{A}\bigcup\{\mathcal{B}_\zeta,\mathcal{C}_\theta~|~\zeta\in (J_f)_{<\lambda}\wedge\theta\in (J_f)_\lambda\wedge\zeta\neq\xi\}.$
 \item $\mathcal{C}_\eta\downarrow_\mathcal{A}\bigcup\{\mathcal{B}_\zeta,\mathcal{C}_\theta~|~\zeta\in (J_f)_{<\lambda}\wedge\theta\in (J_f)_\lambda\wedge\theta\neq\eta\}.$
 \end{itemize}
 We can choose this models due to the existence property (see above) and the finite character.
 Notice that all $\xi\in (J_f)_{<\lambda}$ and $\eta\in (J_f)_{\lambda}$, satisfy $$\mathcal{B}_\xi\cup\mathcal{C}_\eta\downarrow_\mathcal{A}\bigcup\{\mathcal{B}_\zeta,\mathcal{C}_\theta~|~\zeta\in (J_f)_{<\lambda}\wedge\theta\in (J_f)_\lambda\wedge\zeta\neq\xi\wedge\theta\neq\eta\}.$$

Let $F_{\xi\eta}$ be an automorphism of the monster model such that $F_{\xi\eta}\restriction \mathcal{C}:\mathcal{C}\rightarrow\mathcal{C}_\eta$ and $F_{\xi\eta}\restriction \mathcal{B}:\mathcal{B}\rightarrow\mathcal{B}_\xi$ are isomorphisms and $F_{\xi\eta}\restriction\mathcal{A}=id$. Denote the sequence $\mathcal{I}$ by $\{w_\alpha~|~\alpha<\kappa\}$. For all $\eta\in (J_f)_{\lambda}$ and every $\xi<\eta$, let $I_{\xi\eta}=\{b_\alpha~|~\alpha<c_f(\eta)\}$ be an indiscernible sequence over $\mathcal{B}_\xi\cup\mathcal{C}_\eta$
of size $c_f(\eta)$, that is independent over $\mathcal{B}_\xi\cup\mathcal{C}_\eta$, and satisfies:
\begin{itemize}
\item $tp(I_{\xi\eta},\mathcal{B}_\xi\cup\mathcal{C}_\eta)=tp(F_{\xi\eta}(\mathcal{I}\restriction c_f(\eta)),\mathcal{B}_\xi\cup\mathcal{C}_\eta)$.
\item $I_{\xi\eta}\downarrow_{\mathcal{B}_\xi\cup\mathcal{C}_\eta}\bigcup\{\mathcal{B}_\zeta,\mathcal{C}_\theta~|~\zeta\in (J_f)_{<\lambda}\wedge\theta\in (J_f)_\lambda\}\cup\bigcup\{I_{\zeta\theta}~|~\zeta\neq\xi\vee\theta\neq\eta\}$.
\end{itemize}
Therefore, there is an elementary embedding $G:\mathcal{B}_\xi\cup\mathcal{C}_\eta\cup F_{\xi\eta}(\mathcal{I}\restriction c_f(\eta))\rightarrow \mathcal{B}_\xi\cup\mathcal{C}_\eta\cup I_{\xi\eta}$ given by $G\restriction\mathcal{B}_\xi\cup\mathcal{C}_\eta=id$ and $G(F_{\xi\eta}(\mathcal{I}\restriction c_f(\eta)))=I_{\xi\eta}$. So the map $H_{\xi\eta}:\mathcal{B}\cup\mathcal{C}\cup \mathcal{I}\restriction c_f(\eta)\rightarrow \mathcal{B}_\xi\cup\mathcal{C}_\eta\cup I_{\xi\eta}$ given by $H_{\xi\eta}=G\circ (F_{\xi\eta}\restriction dom(H_{\xi\eta}))$ is elementary.
\begin{remark}\label{R.4.4a}
To recap, $\mathcal{B}_\xi$, $\mathcal{C}_\eta$, and $I_{\xi\eta}$ satisfy the following:
\begin{enumerate}
\item [1.]$Av(I_{\xi\eta},\mathcal{B}_\xi\cup\mathcal{C}_\eta)$ is orthogonal to $\mathcal{B}_\xi$ and to $\mathcal{C}_\eta$.
\item [2.]If $(B_i)_{i<3}$ are sets such that:
\begin{itemize}
\item [(a)]$B_0\downarrow_{\mathcal{A}}\mathcal{B}_\xi\cup\mathcal{C}_\eta$.
\item [(b)]$B_1\downarrow_{\mathcal{B}_\xi\cup B_0}B_2\cup\mathcal{C}_\eta$.
\item [(c)]$B_2\downarrow_{\mathcal{C}_\eta\cup B_0}B_1\cup\mathcal{B}_\xi$.
\end{itemize}
Then, $$tp(I_{\xi\eta},\mathcal{B}_\xi\cup\mathcal{C}_\eta) \vdash tp(I_{\xi\eta},\mathcal{B}_\xi\cup\mathcal{C}_\eta\cup_{i<3}B_i).$$
\item [3.]$I_{\xi\eta}\downarrow_{\mathcal{B}_\xi\cup\mathcal{C}_\eta}\bigcup\{\mathcal{B}_\zeta,\mathcal{C}_\theta~|~\zeta\in (J_f)_{<\lambda}\wedge\theta\in (J_f)_\lambda\}\cup\bigcup\{I_{\zeta\theta}~|~\zeta\neq\xi\vee\theta\neq\eta\}$.
\end{enumerate}
\end{remark}
\begin{definition}\label{D.4.5}
Let $\Gamma_f$ be the set $\bigcup\{\mathcal{B}_\xi,\mathcal{C}_\eta,I_{\xi\eta}~|~\xi\in (J_f)_{<\lambda}\wedge\eta\in (J_f)_\lambda \wedge \xi<\eta\}$ and let $\mathcal{A}^f$ be the $a$-primary model over $\Gamma_f$.
Let $\Gamma_f^\alpha$ be the set $\bigcup\{\mathcal{B}_\xi,\mathcal{C}_\eta,I_{\xi\eta}~|~\xi,\eta\in J_f^\alpha \wedge \xi<\eta\}$, recall $J_f^\alpha$ from Remark \ref{R.2.7a}.
\end{definition}
\begin{fact}\label{F.4.6}
If $\alpha$ is such that $\alpha^{\lambda}<f(\alpha)$, $sup(\{c_f(\eta)\}_{\eta\in J_f^\alpha})<\alpha$, then $|\Gamma_f^{\alpha+1}|= f(\alpha)$.
\end{fact}
\begin{proof}
Recall that for all $\xi\in (J^\alpha_f)_{<\lambda}$ and $\eta\in (J^\alpha_f)_\lambda $, $\mathcal{B}_\xi$ and $\mathcal{C}_\eta$ have cardinality $2^\omega$, and $I_{\xi\eta}$ has cardinality $c_f(\eta)$.
Since $\Gamma_f^\alpha=\cup\{\mathcal{B}_\xi,\mathcal{C}_\eta,I_{\xi\eta}~|~\xi\in (J^\alpha_f)_{<\lambda}\wedge\eta\in (J^\alpha_f)_\lambda \wedge \xi<\eta\}$, we know that $|\Gamma_f^{\alpha+1}|\leq |J^{\alpha+1}_f|\cdot sup(\{c_f(\eta)\}_{\eta\in (J_f^{\alpha+1})_{\lambda}})$. Since $|J_f^{\alpha+1}|\leq \alpha^{\lambda}<f(\alpha)$ and $sup(\{c_f(\eta)\}_{\eta\in J_f^\alpha})<\alpha<f(\alpha)$, we get $|\Gamma_f^{\alpha+1}|\leq max(f(\alpha), sup(\{c_f(\eta)\}_{\eta\in J_f^{\alpha+1}\backslash J_f^\alpha}))$. But every $\eta\in J_f^{\alpha+1}\backslash J_f^\alpha$ with domain $\lambda$ has $rang(\eta_1)=\lambda$ and $f(\alpha)=c_f(\eta)$, otherwise $sup(rang(\eta_5))<\alpha$ and $\eta\in J_f^\alpha$. We conclude $|\Gamma_f^{\alpha+1}|= f(\alpha)$.
\end{proof}

\begin{lemma}\label{L.4.7}
For every $\xi\in (J_f)_{<\lambda}$, $\eta\in (J_f)_\lambda$, $\xi<\eta$, let $p_{\xi\eta}$ be the type $Av(I_{\xi\eta}\restriction \omega,I_{\xi\eta}\restriction\omega\cup\mathcal{B}_\xi\cup\mathcal{C}_\eta)$. If $c_f(\eta)>\omega$, then $dim(p_{\xi\eta},\mathcal{A}^f)=c_f(\eta)$.
\end{lemma}
\begin{proof}
Denote by $S$ the set $I_{\xi\eta}\restriction\omega\cup\mathcal{B}_\xi\cup\mathcal{C}_\eta$, so $p_{\xi\eta}=Av(I_{\xi\eta}\restriction\omega,S)$.

Suppose, towards a contradiction, that $dim(p_{\xi\eta},\mathcal{A}^f)\neq c_f(\eta)$. Since $I_{\xi\eta}\subset \mathcal{A}^f$, then $dim(p_{\xi\eta},\mathcal{A}^f)> c_f(\eta)$. Therefore, there is an independent sequence $I=\{a_i~|~i<c_f(\eta)^+\}$ over $S$ such that $I\subset \mathcal{A}^f$ and $\forall a\in I$, $a\models p_{\xi\eta}$.

By induction on $\alpha$, it can be proved that $I_{\xi\eta}\restriction\omega\cup \{a_i~|~i\leq \alpha\}$ is indiscernible over $\mathcal{B}_\xi\cup \mathcal{C}_\eta$. Therefore $I_{\xi\eta}\restriction\omega\cup I$ is indiscernible over $\mathcal{B}_\xi\cup \mathcal{C}_\eta$.
In particular $I_{\xi\eta}\restriction\omega\cup I$ is indiscernible, and $I_{\xi\eta}$ is equivalent to $I$.

By some forking calculus manipulation and Remark \ref{R.4.4a}, it is easy to prove that 
$tp(I_{\xi\eta},\mathcal{B}_\xi\cup\mathcal{C}_\eta)\vdash tp(I_{\xi\eta},\Gamma_f\backslash I_{\xi\eta})$ and $I_{\xi\eta}$ is indiscernible over $\Gamma_f\backslash I_{\xi\eta}$.

We know that $tp(I,\mathcal{B}_\xi\cup\mathcal{C}_\eta)=tp(I_{\xi\eta},\mathcal{B}_\xi\cup\mathcal{C}_\eta)$, therefore 
$tp(I,\mathcal{B}_\xi\cup\mathcal{C}_\eta)\vdash tp(I_{\xi\eta},\Gamma_f\backslash I_{\xi\eta})$. We conclude that $tp(I,\mathcal{B}_\xi\cup\mathcal{C}_\eta)\vdash tp(I,\Gamma_f\backslash I_{\xi\eta})$ and since $I$ is indiscernible over $\mathcal{B}_\xi\cup\mathcal{C}_\eta$, then $I$ is indiscernible over $\Gamma_f\backslash I_{\xi\eta}$.

There are $I',I^*\subseteq I$ such that $|I'|=c_f(\eta)^+$ and $I'\downarrow_{(\Gamma_f\backslash I_{\xi\eta})\cup I^*}I_{\xi\eta}$. 

In particular $I'$ is indiscernible over $\Gamma_f\cup I^*$, and $I'$ is indiscernible over $\Gamma_f$.

Let $J\subset \mathcal{A}^f$ be a maximal indiscernible set over $\Gamma_f$ such that $I'\subseteq J$. By Lemma \ref{L.4.2} $|J|+\kappa(T)=dim(J,\Gamma_f,\mathcal{A}^f)+\kappa(T)$. Since $T$ is superstable, $\kappa(T)=\omega<c_f(\eta)^+\leq |J|$ and we conclude that  $\kappa(T)<dim(J,\Gamma_f,\mathcal{A}^f)+\kappa(T)$. Therefore $\kappa(T)<dim(J,\Gamma_f,\mathcal{A}^f)$ and by Lemma \ref{L.4.2} the dimension is true, $dim(J,\Gamma_f,\mathcal{A}^f)=|J|$. So $dim(J,\Gamma_f,\mathcal{A}^f)>\omega$ a contradiction with Theorem \ref{T.4.3}.
\end{proof}

\begin{lemma}\label{lemma_no_sequence}
For every $\xi\in (J_f)_{<\lambda}$, $\eta\in (J_f)_\lambda$, $\xi<\eta$ and $\alpha$ such that $c(\eta)=f(\alpha)$, where $f(\alpha)$ is a cardinal, there is no independent indiscernible sequence $(c_i)_{i<f(\alpha)^+}$ over $\mathcal{B}_\xi\cup\mathcal{C}_\eta$ with $c_0\in I_{\xi\eta}$. 
\end{lemma}
\begin{proof}

Denote by $J$ the sequence $(c_i)_{i<f(\alpha)^+}$, since $T$ is superstable, there is $J'\subseteq J$ of power $f(\alpha)^+$ such that $c_0\notin J'$ and satisfies $J'\downarrow_{J\restriction\omega\cup\mathcal{B}_\xi\cup\mathcal{C}_\eta}I_{\xi\eta}$.
Since $J$ is an independent sequence over $\mathcal{B}_\xi\cup\mathcal{C}_\eta$, then $J'\downarrow_{\mathcal{B}_\xi\cup\mathcal{C}_\eta}J\restriction\omega\cup I_{\xi\eta}.$ 
Let us denote by $Q$ the set $\mathcal{B}_\xi\cup\mathcal{C}_\eta\cup (I_{\xi\eta}\restriction\omega)\backslash \{c_0\}$,
so $J'\downarrow_{Q}I_{\xi\eta}$. Since $Av(I_{\xi\eta},Q)$ is stationary and $I_{\xi\eta}$ is independent over $\mathcal{B}_\xi\cup\mathcal{C}_\eta$, we conclude that $I'=\{c_0\}\cup(I_{\xi\eta}\backslash (I_{\xi\eta}\restriction\omega))$ is indiscernible over $J'\cup Q$. Especially $I'$ is indiscernible over $\mathcal{B}_\xi\cup\mathcal{C}_\eta\cup J'$. On the other hand $J'\downarrow_{\mathcal{B}_\xi\cup\mathcal{C}_\eta}J\restriction\omega\cup I_{\xi\eta}$ implies that $J'\downarrow_{\mathcal{B}_\xi\cup\mathcal{C}_\eta}I_{\xi\eta}$, and since $I_{\xi\eta}$ is independent over $\mathcal{B}_\xi\cup\mathcal{C}_\eta$, we conclude that $I_{\xi\eta}$ is independent over $\mathcal{B}_\xi\cup\mathcal{C}_\eta\cup J'$. In particular $I'$ is independent over $\mathcal{B}_\xi\cup\mathcal{C}_\eta\cup J'$.
\begin{claim}\label{inductionseq}
$J'\cup I'$ is indiscernible over $\mathcal{B}_\xi\cup\mathcal{C}_\eta$
\end{claim}
\textit{Proof of Claim \ref{inductionseq}.}
We will proceed by induction. Let us denote by $\{d_i\mid i<f(\alpha)\}$ the sequence $I'$. 
Since $c_0\in I'\cap J$, $c_0\models Av(J',\mathcal{B}_\xi\cup\mathcal{C}_\eta\cup J')$, and $I'$ is indiscernible over $J'\cup Q$, then for every $i<f(\alpha)$,
\begin{equation*}
d_i\models Av(J',\mathcal{B}_\xi\cup\mathcal{C}_\eta\cup J').
\end{equation*}
Suppose $j$ is such that for all $n<j$ the sequence $J'\cup \{d_i\mid i\leq n\}$ is indiscernible over $\mathcal{B}_\xi\cup\mathcal{C}_\eta$, then $J'\cup \{d_i\mid i<j\}$ is indiscernible over $\mathcal{B}_\xi\cup\mathcal{C}_\eta$, therefore $Av(J'\cup \{d_i\mid i<j\},\mathcal{B}_\xi\cup\mathcal{C}_\eta\cup J'\cup \{d_i\mid i<j\})=Av(J',\mathcal{B}_\xi\cup\mathcal{C}_\eta\cup J'\cup \{d_i\mid i<j\})$ and it does not fork over $\mathcal{B}_\xi\cup\mathcal{C}_\eta\cup J'$. On the other hand we know that $Av(J',\mathcal{B}_\xi\cup\mathcal{C}_\eta\cup J')$ is stationary,
$d_j\downarrow_{\mathcal{B}_\xi\cup\mathcal{C}_\eta\cup J'}\{d_i\mid i<j\}$ and $d_j\models Av(J',\mathcal{B}_\xi\cup\mathcal{C}_\eta\cup J')$, we conclude that $tp(d_j,\mathcal{B}_\xi\cup\mathcal{C}_\eta\cup J'\cup \{d_i\mid i<j\}))=Av(J'\cup \{d_i\mid i<j\},\mathcal{B}_\xi\cup\mathcal{C}_\eta\cup J'\cup \{d_i\mid i<j\})$. Therefore $J'\cup \{d_i\mid i\leq j\}$ is indiscernible over $\mathcal{B}_\xi\cup\mathcal{C}_\eta$. We conclude that $J'\cup I'$ is indiscernible over $\mathcal{B}_\xi\cup\mathcal{C}_\eta$. This finishes the proof of Claim \ref{inductionseq}.

So $J'$ is equivalent to $I_{\xi\eta}$ and for all $d\in J'$, $d\models Av(I_{\xi\eta}\restriction \omega,I_{\xi\eta}\restriction \omega\cup \mathcal{B}_\xi\cup\mathcal{C}_\eta)$. Since $J'$ is independent over $\mathcal{B}_\xi\cup\mathcal{C}_\eta$ and $J'\downarrow_{\mathcal{B}_\xi\cup\mathcal{C}_\eta}I_{\xi\eta}$, we conclude that $J'$ is independent over $I_{\xi\eta}\restriction \omega\cup \mathcal{B}_\xi\cup\mathcal{C}_\eta$, thus $dim(p_{\xi\eta},\mathcal{A}^f)\ge f(\alpha)^+$, but this contradicts Lemma \ref{L.4.7}.
\end{proof}

One of the key lemmas for the following section
 is Lemma \ref{L.4.10} (below). 
Let us define the nice subsets of $\Gamma_f$. These subsets have a couple of properties, that will be useful when we study the model $\mathcal{A}^f$.
\begin{definition}\label{D.4.9}
We say $\emptyset\neq X\subseteq \Gamma_f$ is nice if the following holds.
\begin{enumerate}
\item [1.]If $X\cap I_{\xi\eta}\neq \emptyset$, then $\mathcal{B}_\xi,\mathcal{C}_\eta\subset X$.
\item [2.]If $\mathcal{B}_\xi\cap X\neq \emptyset$, then $\mathcal{B}_\xi\subset X$.
\item [3.]If $\mathcal{C}_\eta\cap X\neq \emptyset$, then $\mathcal{C}_\eta\subset X$.
\item [4.]If $\xi<\eta$ and $\mathcal{B}_\xi,\mathcal{C}_\eta\subset X$, then $X\cap I_{\xi\eta}$ is infinite.
\end{enumerate}
\end{definition}
The argument for the following Lemma is a variation of the argument used in \cite{HS99} in the fourth section. 
\begin{lemma}\label{L.4.10}
Let $Z$ be a nice subset of $\Gamma_f$ and $d\in\Gamma_f\backslash Z$. Then for all $B$ finite subset of $Z$ there is $f\in Saut(\mathcal{M},B)$ such that $f(d)\in Z$.
\end{lemma}

Suppose $X$ and $A$ are nice subsets of $\Gamma_f$. If $\xi$ and $\eta$ are such that $\mathcal{B}_\xi\cup\mathcal{C}_\eta\subseteq A$ and $I_{\xi\eta}\cap X\subseteq A$, then we say that $A$ is $X$-nice for $(\xi,\eta)$. 
\begin{lemma}\label{L.4.11}
Suppose $Z\subseteq \Gamma_f$ is nice and $B$ is $a$-constructable over $Z$. If $X\subseteq \Gamma_f$ is a nice subset such that $Z\cup X$ is nice, then  $B\cup X$ is $a$-constructible over $Z\cup X$.
\end{lemma}
\begin{proof}
Let $(Z,(a_i,B_i)_{i<\gamma})$ be an $a$-construction for $B$ over $Z$. Let $(\mathcal{D}_i)_{i<\delta}$ be an enumeration of $\{\mathcal{B}_\xi,\mathcal{C}_\eta,I_{\xi\eta}\cap X~|~\xi<\eta\wedge \mathcal{B}_\xi\cup\mathcal{C}_\eta\subseteq Z\cup X\}$ such that $\mathcal{B}_\xi$ and $\mathcal{C}_\eta$ are before $I_{\xi\eta}$ in the enumeration. Let $Z^j$ be the minimal nice subset of $Z\cup X$ that contains $Z\cup\bigcup_{i\leq j}\mathcal{D}_i$, and it is $X$-nice for every $(x,y)$ that satisfies: either $\mathcal{B}_x\subseteq\bigcup_{i\leq j}\mathcal{D}_i\backslash Z$ or $\mathcal{C}_y\subseteq\bigcup_{i\leq j}\mathcal{D}_i\backslash Z$. 
\begin{claim}\label{claim.4.11}
$(Z^j,(a_i,B_i)_{i<\gamma})$ is an $a$-construction for $B\cup Z^j$ over $Z^j$, for every $j<\delta$.
\end{claim}
\textit{Proof of Claim \ref{claim.4.11}.}
Suppose, towards a contradiction, that $\alpha$ is the minimal ordinal such that $(Z^\alpha,(a_i,B_i)_{i<\gamma})$ is not an $a$-construction for $B\cup Z^\alpha$ over $Z^\alpha$.
By the minimality of $\alpha$, $(Z^\beta,(a_i,B_i)_{i<\gamma})$ is an $a$-construction for $B\cup Z^\beta$ over $Z^\beta$, for every $\beta<\alpha$. Therefore for every $\beta<\alpha$ and $i<\gamma$, $(tp(a_i,Z_i^\beta),B_i)\in F_\omega^a$ where $Z_i^\beta=Z^\beta\cup\bigcup_{j<i}a_j$. So $(tp(a_i,\cup_{\beta<\alpha}Z_i^\beta),B_i)\in F_\omega^a$ for every $i<\gamma$, we conclude that $\alpha$ is not a limit cardinal. Let us denote by $Z'$ the set $Z^\beta$, for $\beta$ the predecessor of $\alpha$.

The proof is divided in the following cases:
\begin{enumerate}
\item [1.]$\mathcal{D}_\alpha=\mathcal{C}_\eta$ for some $\mathcal{C}_\eta\subseteq X\cup Z$.
\item [2.]$\mathcal{D}_\alpha=\mathcal{B}_\xi$ for some $\mathcal{B}_\xi\subseteq X\cup Z$.
\item [3.]$\mathcal{D}_\alpha=I_{\xi\eta}\cap X$, for some $\mathcal{B}_\xi\cup\mathcal{C}_\eta\subseteq X\cup Z$.
\end{enumerate}
The case 2 is similar to the case 1, we will show only the cases 1 and 3.

Case 1: 
Since $(Z^\alpha,(a_i,B_i)_{i<\gamma})$ is not an $a$-construction over $Z^\alpha$, then by the minimality of $Z^\alpha$, $\mathcal{C}_\eta\not \subseteq Z'$. Therefore, $I_{\xi\eta}\cap Z'=\emptyset$ for every $\xi<\eta$. Since $X\cup Z$ is nice, then we know that for all $\mathcal{B}_\xi\subseteq Z'$ that satisfies $\xi<\eta$, it holds that $\mathcal{B}_\xi\subseteq X$.
Let $n$ be the least ordinal such that $(Z'\cup \mathcal{C}_\eta\cup \bigcup \{I_{\xi\eta}\cap X~|~\xi<\eta\wedge \mathcal{B}_\xi\subseteq Z'\},(a_i,B_i)_{i\leq n})$ is not an a-construction over $Z'\cup \mathcal{C}_\eta\cup \bigcup \{I_{\xi\eta}\cap X~|~\xi<\eta\wedge \mathcal{B}_\xi\subseteq Z'\}$, since $a$-isolation is the $F_\omega^a$-isolation, then $B_n$ is finite and we can assume $n<\omega$.

Denote by $D$ the set $\mathcal{C}_\eta\cup \bigcup \{I_{\xi\eta}\cap X~|~\xi<\eta\wedge \mathcal{B}_\xi\subseteq Z'\}$.
Since $(Z'\cup D,(a_i,B_i)_{i< n})$ is an $a$-construction over $Z'$, then $C=\bigcup_{i< n} B_i\cap (Z'\cup D)$ is such that $stp(a_0^\frown\cdots ^\frown a_{n-1},C)\vdash tp(a_0^\frown\cdots ^\frown a_{n-1}, Z'\cup D)$. Notice that $C$ is a subset of $Z'$.
On the other hand, there is $b$ such that $stp(b,B_n)=stp(a_n,B_n)$, and $tp(b,Z'\cup \bigcup\{a_i~|~i<n\}\cup D)\neq tp(a_n,Z'\cup \bigcup\{a_i~|~i<n\}\cup D)$. So there are tuples $d\in D\backslash \mathcal{A}$ and $e\in Z'\cup \bigcup\{a_i~|~i<n\}$ that satisfy $tp(b,e\cup d)\neq tp(a_n,e\cup d)$. Denote by $W$ the set $C\cup((B_n\cup e)\cap Z')$, by Lemma \ref{L.4.10} we know that there is $g\in Saut(\mathcal{M},W)$ such that $g(d)\in Z'$. We know that, $stp(a_0^\frown\cdots ^\frown a_{n-1},C)\vdash tp(a_0^\frown\cdots ^\frown a_{n-1}, Z'\cup D)$, so $a_0^\frown\cdots ^\frown a_{n-1}\downarrow_C Z'\cup D$. We conclude that $$a_0^\frown\cdots ^\frown a_{n-1}\downarrow_W d$$ and $$a_0^\frown\cdots ^\frown a_{n-1}\downarrow_W g(d).$$ Therefore $stp(d,C\cup B_n\cup e)=stp(g(d), C\cup B_n\cup e)$ and there is $f\in Saut(\mathcal{M},C\cup B_n\cup e)$ that satisfies $f(d)=g(d)$.

Since $tp(b,e\cup d)\neq tp(a_n, e\cup d)$ and $stp(b,B_n)=stp(a_n,B_n)$ hold, then we have that $tp(f(b),e\cup f(d))\neq tp(f(a_n), e\cup f(d))$, and the strong types of $a_n,b,f(a_n)$ and $f(b)$ over $B_n$ are the same strong type. Since $(Z',(a_i,B_i)_{i<\gamma})$ is an $a$-construction, then by the $a$-isolation we know that $stp(a,B_n)\vdash tp(a_n,Z'\cup \bigcup\{a_i~|~i<n\})$, on the other hand $stp(a_n,B_n)=stp(f(a_n),B_n)=stp(f(b),B_n)$, so $tp(f(a_n),Z'\cup \bigcup\{a_i~|~i<n\})=tp(f(b),Z'\cup \bigcup\{a_i~|~i<n\})$. In particular $e,f(d)\in Z'$, so $tp(f(b),e\cup f(d))=tp(f(a_n),e\cup f(d))$, a contradiction.

Case 3: 
By the way $(\mathcal{D}_i)_{i<\delta}$ was defined, we know that $\mathcal{B}_\xi$ and $\mathcal{C}_\eta$ are before $I_{\xi\eta}\cap X$ in the enumeration, so $\mathcal{B}_\xi\cup\mathcal{C}_\xi\subseteq Z'$. We have the following possibilities, either $\mathcal{B}_\xi\not\subseteq Z$, or $\mathcal{C}_\eta\not\subseteq Z$, or $\mathcal{B}_\xi,\mathcal{C}_\eta\subseteq Z$. In the first two cases, by the way $Z'$ was defined, we know that $Z'$ is $X$-nice for $(\xi, \eta)$, so $I_{\xi\eta}\cap X\subset Z'$. Therefore, $Z'=Z^\alpha$ and $(Z',(a_i,B_i)_{i<\gamma})$ is an $a$-construction for $B\cup Z^\alpha$ over $Z^\alpha$, a contradiction. Therefore, we need to show only the case when $\mathcal{B}_\xi,\mathcal{C}_\eta\subset Z$. Since $(Z^\alpha,(a_i,B_i)_{i<\gamma})$ is not an $a$-construction over $Z^\alpha$, then $I_{\xi\eta}\cap X\not \subseteq Z'$.

Let $n$ be the least ordinal such that $(Z'\cup (I_{\xi\eta}\cap X),(a_i,B_i)_{i\leq n})$ is not an a-construction over $Z'\cup (I_{\xi\eta}\cap X)$, since $a$-isolation is the $F_\omega^a$-isolation, then $B_n$ is finite and we can assume $n<\omega$.
Since $(Z'\cup (I_{\xi\eta}\cap X),(a_i,B_i)_{i< n})$ is an $a$-construction over $Z'\cup (I_{\xi\eta}\cap X)$, then $C=\bigcup_{i< n} B_i\cap (Z'\cup (I_{\xi\eta}\cap X))$ is such that $stp(a_0^\frown\cdots ^\frown a_{n-1},C)\vdash tp(a_0^\frown\cdots ^\frown a_{n-1}, Z'\cup( I_{\xi\eta}\cap X))$. Notice that $C$ is a subset of $Z'$.
On the other hand, there is $b$ such that $stp(b,B_n)=stp(a_n,B_n)$, and $tp(b,Z'\cup \bigcup\{a_i~|~i<n\}\cup (I_{\xi\eta}\cap X))\neq tp(a_n,Z'\cup \bigcup\{a_i~|~i<n\}\cup (I_{\xi\eta}\cap X))$. Since $Z'$ is nice, then there is an infinite $I'_{\xi\eta}\subset I_{\xi\eta}\cap X$ contained in $Z'$.
Therefore, there are tuples $d\in (I_{\xi\eta}\cap X)\backslash I'_{\xi\eta}$ and $e\in Z'\cup \bigcup\{a_i~|~i<n\}$ that satisfy $tp(b,e\cup d)\neq tp(a_n, e\cup d)$. Denote by $W$ the set $C\cup ((B_n\cup e)\cap Z')$, by Lemma \ref{L.4.10} we know that there is $g\in Saut(\mathcal{M},W)$ such that $g(d)\in Z'$.
Since $stp(a_0^\frown\cdots ^\frown a_{n-1},C)\vdash tp(a_0^\frown\cdots ^\frown a_{n-1}, Z'\cup (I_{\xi\eta}\cap X))$, then $a_0^\frown\cdots ^\frown a_{n-1}\downarrow_C Z'\cup (I_{\xi\eta}\cap X)$. Therefore $$a_0^\frown\cdots ^\frown a_{n-1}\downarrow_W d$$ and $$a_0^\frown\cdots ^\frown a_{n-1}\downarrow_W g(d).$$
So, $stp(d,C\cup B_n\cup e)=stp(g(d), C\cup B_n\cup e)$ and there is $f\in Saut(\mathcal{M},C\cup B_n\cup e)$ that satisfies $f(d)=g(d)$. From here, we finish as in case 1.  

This finishes the proof of Claim \ref{claim.4.11}.


Since for every $\beta<\delta$ and $i<\gamma$, $(tp(a_i,Z_i^\beta),B_i)\in F_\omega^a$ where $Z_i^\beta=Z^\beta\cup\bigcup_{j<i}a_j$, then $(tp(a_i,\cup_{\beta<\delta}Z_i^\beta),B_i)\in F_\omega^a$ and $(Z\cup X,(a_i,B_i)_{i<\gamma})$ is an $a$-construction for $B\cup X$ over $Z\cup X$.
\end{proof}
\begin{fact}\label{F.4.12}
If $Z\subseteq \Gamma_f$ is nice, then for every $\alpha<\kappa$ the following holds $$Z\downarrow_{Z\cap \Gamma_f^\alpha}\Gamma_f^\alpha.$$ 
\end{fact}

\begin{corollary}\label{C.4.13}
If $Z\subseteq \Gamma_f$ is nice, then for every nice set $\Gamma\subseteq \Gamma_f$ the following holds $$Z\downarrow_{Z\cap \Gamma}\Gamma.$$ 
\end{corollary}
Now, we have all the tools needed to prove the main result on $\mathcal{A}^f$. 

\section{Main result on $\mathcal{A}^f$}
This section is devoted to prove, for a certain kind of functions, that the models $\mathcal{A}^f$ and $\mathcal{A}^g$ are isomorphic if and only if $J_f$ and $J_g$ are isomorphic coloured trees.
\begin{theorem}\label{T.4.14}
Assume $f,g$ are functions from $\kappa$ to $Card\cap\kappa\backslash \lambda$ such that $f(\alpha),g(\alpha)>\alpha^{++}$ and $f(\alpha),g(\alpha)>\alpha^\lambda$. Then $\mathcal{A}^f$ and $\mathcal{A}^g$ are isomorphic if and only if $f$ and $g$ are $E_{\lambda\text{-club}}^\kappa$ equivalent.
\end{theorem}

\begin{lemma}
Assume $f,g$ are functions from $\kappa$ to $\kappa$. 
If $f$ and $g$ are $E_{\lambda\text{-club}}^\kappa$ equivalent, then $\mathcal{A}^f$ and $\mathcal{A}^g$ are isomorphic.
\end{lemma}
\begin{proof}
Assume $f$ and $g$ are $E_{\lambda\text{-club}}^\kappa$ equivalent. By Lemma \ref{L.2.7} $J_f$ and $J_g$ are isomorphic coloured trees, let $G:J_f\rightarrow J_g$ be an isomorphism. Define $\mathcal{H}_{\xi\eta}:\mathcal{B}_\xi\cup\mathcal{C}_\eta\cup I_{\xi\eta}\rightarrow\mathcal{B}_{G(\xi)}\cup\mathcal{C}_{G(\eta)}\cup I_{G(\xi)G(\eta)}$ by $\mathcal{H}_{\xi\eta}=H_{G(\xi)G(\eta)}\circ H^{-1}_{\xi\eta}$ (where $H_{rp}$ is the elementary embedding used in the construction of $I_{rp}$), clearly $\mathcal{H}_{\xi\eta}$ is elementary.
It is easy to check that the map $$\mathcal{H}=\bigcup_{\eta\in (J_f)_\lambda}\bigcup_{\xi\in (J_f)_{<\lambda}, \xi<\eta}\mathcal{H}_{\xi\eta}$$ is elementary.
Let $\bar{\mathcal{H}}$ be an automorphism that extends $\mathcal{H}$, then $\bar{\mathcal{H}}(\mathcal{A}^f)$ is $a$-primary over $\Gamma_g$. Therefore $\bar{\mathcal{H}}(\mathcal{A}^f)$ and $\mathcal{A}^g$ are isomorphic, we conclude that $\mathcal{A}^f$ and $\mathcal{A}^g$ are isomorphic.

\end{proof}
\begin{lemma}
Assume $f,g$ are functions from $\kappa$ to $Card\cap\kappa\backslash \lambda$ such that $f(\alpha),g(\alpha)>\alpha^{++}$ and $f(\alpha),g(\alpha)>\alpha^\lambda$. If $\mathcal{A}^f$ and $\mathcal{A}^g$ are isomorphic, then $f$ and $g$ are $E_{\lambda\text{-club}}^\kappa$ equivalent.
\end{lemma}

{\bf Comments: }This lemma has a long proof and it is divided in five claim. In the proof we proceed by contradiction, we assume that $f$ and $g$ are not $E_{\lambda\text{-club}}^\kappa$ equivalent and there is an isomorphism $\Pi:\mathcal{A}^f\rightarrow \mathcal{A}^g$. Then we construct an $a$-primary model $F$, and find $\xi<\eta$ and $a\in I_{\xi\eta}$ such that $$\Pi(a)\downarrow_{\Pi(\mathcal{B}_\xi\cup\mathcal{C}_\eta)}F.$$ By using $a$, we will construct a independent indiscernible sequence $(b_i)_{i<f(\alpha)^+}$ over $\Pi(\mathcal{B}_\xi\cup\mathcal{C}_\eta)$ in $\mathcal{A}^g$. Finally, we use $\Pi$ and $(b_i)_{i<f(\alpha)^+}$ to construct a sequence $(c_i)_{i<f(\alpha)^+}$ indiscernible and independent over $\mathcal{B}_\xi\cup\mathcal{C}_\eta$ with $c_0\in I_{\xi\eta}$, which is a contradiction with Lemma \ref{lemma_no_sequence}.

\begin{proof}

Let us assume, towards a contradiction, that $f$ and $g$ are not $E_{\lambda\text{-club}}^\kappa$ equivalent and there is an isomorphism $\Pi:\mathcal{A}^f\rightarrow \mathcal{A}^g$. Without loss of generality, we can assume that $\{\alpha~|~f(\alpha)>g(\alpha)\wedge cf(\alpha)=\lambda\}$ is stationary.
Let $(\Gamma_f,(a_i^f,B_i^f)_{i<\gamma})$ be an $a$-construction of $\mathcal{A}^f$ over $\Gamma_f$. For every $\alpha$ define $\mathcal{A}_f^\alpha=\Gamma^\alpha_f\cup\bigcup\{a_i^f~|~i<\alpha\}$, clearly $\mathcal{A}_f^\alpha$ is not necessary a model.

We say that $\alpha<\kappa$ is $f$-good if $(\Gamma^\alpha_f,(a_i^f,B_i^f)_{i<\alpha})$ is an $a$-construction over $\Gamma^\alpha_f$, $\mathcal{A}_f^\alpha$ is an $a$-primary model over $\Gamma_f^\alpha$, and $\alpha$ is a cardinal. Notice that there are club many $f$-good cardinals.
We say that $\alpha$ is very good if, $\alpha$ is $f$-good, $f(\alpha)>g(\alpha)>\alpha^{++}$ and $\Pi(\mathcal{A}_f^\alpha)=\mathcal{A}_g^\alpha$. 
\begin{claim}\label{claim1final}
There is a very good cardinal $\alpha$ with cofinality $\lambda$ and $\eta\in J_f^{\alpha+1}$  such that 
\begin{itemize}
\item $\alpha^{\lambda}<g(\alpha)$,
\item $sup(\{c_g(p)\}_{p\in J_g^\alpha})<\alpha$,
\item there are cofinally many very good cardinals $\beta<\alpha$,
\item $\bigcup rang(\eta_1)=\lambda$ and $\bigcup rang(\eta_5)=\alpha$.
\end{itemize}
\end{claim}
\textit{Proof of Claim \ref{claim1final}.}
Since there are club many $\alpha$'s satisfying $\Pi(\mathcal{A}_f^\alpha)=\mathcal{A}_g^\alpha$ and stationary many $\alpha$'s with cofinality $\lambda$ such that $f(\alpha)>g(\alpha)$, there are stationary many very good cardinals.
Since there are club many $\alpha$'s satisfying $sup(\{c_g(p)\}_{p\in J_g^\alpha})<\alpha$, by Remark \ref{R.2.8} we can choose $\alpha$ a very good cardinal with cofinality $\lambda$ and $\eta\in J_f$ with the desire properties. This finishes the proof of Claim \ref{claim1final}.

Notice that by Definition \ref{D.2.6} item 10, $c_f(\eta)=f(\alpha)$.
Let us choose $X\subseteq \Gamma_g$ and $Y\subseteq \gamma$ such that:
\begin{itemize}
\item $Y$ has power $2^\omega$ and is closed (i.e. for all $i\in Y$, $B_i^g\subseteq \Gamma_g\cup\bigcup_{j\in Y}a^g_j$).
\item $X$ has power $2^\omega$ and is nice.
\item $D=X\cup\bigcup\{a_i^g~|~i\in Y\}$ is the $a$-primary model over $X$.
\item $D^\alpha=(X\cap \Gamma_g^\alpha)\cup\bigcup\{a_i^g~|~i\in Y\wedge i<\alpha\}$ is the $a$-primary model over $X\cap \Gamma_g^\alpha$.
\item $\Pi(\mathcal{C}_\eta)\subseteq D$ and $\Pi(\mathcal{A})\subseteq D^\alpha$.
\item If $\xi\in (J_g)_{<\lambda}$ is such that $\mathcal{B}_\xi\subseteq X$, then for all $\zeta<\xi$, $\mathcal{B}_\zeta\subseteq X$.
\item If $\theta\in (J_g)_{\lambda}\backslash J_g^{\alpha+1}$ is such that $\mathcal{C}_\theta\subseteq X$, then for all $\zeta\in J_g^\alpha$, $\zeta<\theta$ implies that $\mathcal{B}_\zeta\subseteq X$.
\end{itemize}

Let $E$ be an $a$-primary model over $\Gamma_g^{\alpha+1}\cup \mathcal{A}_g^\alpha\cup D$, $F$ be an $a$-primary model 
over $E \cup \bigcup \{\mathcal{B}_\xi, I_{\xi\theta}~|~\xi<\theta\wedge \mathcal{C}_\theta\subseteq X\backslash\Gamma_g^{\alpha+1}\}$, and $G$ be an $a$-primary model over $\Gamma_g\cup F$.
\begin{claim}\label{claim2final}
$E$ is $a$-constructable over $\Gamma_g^{\alpha+1}\cup X$, $F$ is $a$-constructable over $\Gamma_g^{\alpha+1}\cup X\cup \bigcup \{\mathcal{B}_\xi, I_{\xi\theta}~|~\xi<\theta\wedge \mathcal{C}_\theta\subseteq X\backslash\Gamma_g^{\alpha+1}\}$, and without loss of generality, we can assume $G=\mathcal{A}^g$.
\end{claim}
\textit{Proof of Claim \ref{claim2final}.}
Notice that since $D=X\cup\bigcup\{a_i^g~|~i\in Y\}$ is an $a$-construction over $X$, then for all $i\in Y$, $B_i^g\subseteq X\cup\bigcup_{j\in Y}a^g_j$ holds.
 By the definition of $\mathcal{A}^g$, we know that $stp(a_i^g,B_i^g)\vdash tp(a_i^g,\Gamma_g\cup \bigcup\{a_j^g~|~j<i\})$. Since $B_i^g\subseteq X\cup\bigcup\{a_j^g~|~j<i\wedge j\in Y\}$ holds for every $i\in Y$, then $stp(a_i^g,B_i^g)\vdash tp(a_i^g,X\cup \Gamma^\alpha_g\cup \bigcup\{a_j^g~|~j<\alpha\}\cup \bigcup\{a_j^g~|~j<i\wedge j\in Y\})$ holds for all $i\in Y\backslash \alpha$. We conclude that $D\cup \mathcal{A}_g^\alpha$ is $a$-constructable over $X\cup \mathcal{A}_g^\alpha$. Notice that $X\cup\Gamma_g^\alpha$ is nice, so by Lemma \ref{L.4.11} $X\cup \mathcal{A}_g^\alpha$ is $a$-constructable over $X\cup\Gamma_g^\alpha$. We conclude by Lemma \ref{L.4.11} that $E$ is $a$-constructable over $\Gamma_g^{\alpha+1}\cup X$. Notice that $\Gamma_g^{\alpha+1}\cup X\cup \bigcup \{\mathcal{B}_\xi, I_{\xi\theta}~|~\xi<\theta\wedge \mathcal{C}_\theta\subseteq X\backslash\Gamma_g^{\alpha+1}\}$ is nice and by Lemma \ref{L.4.11} we conclude that $F$ is $a$-constructable over $\Gamma_g^{\alpha+1}\cup X\cup \bigcup \{\mathcal{B}_\xi, I_{\xi\theta}~|~\xi<\theta\wedge \mathcal{C}_\theta\subseteq X\backslash\Gamma_g^{\alpha+1}\}$. Since $F$ is $a$-constructable over $\Gamma_g^{\alpha+1}\cup X \cup \bigcup \{\mathcal{B}_\xi, I_{\xi\theta}~|~\xi<\theta\wedge \mathcal{C}_\theta\subseteq X\backslash\Gamma_g^{\alpha+1}\}$, then by Lemma \ref{L.4.11} $G$ is $a$-primary over $\Gamma_g^{\alpha+1}\cup X\cup\bigcup \{\mathcal{B}_\xi, I_{\xi\theta}~|~\xi<\theta\wedge \mathcal{C}_\theta\subseteq X\backslash\Gamma_g^{\alpha+1}\}\cup \Gamma_g$. Therefore, without loss of generality, we can assume $G=\mathcal{A}^g$. This finishes the proof of Claim \ref{claim2final}.
 
 \begin{claim}\label{claim3final}
 There are $\xi<\eta$ and $a\in I_{\xi\eta}\backslash (I_{\xi\eta}\restriction \omega)$ such that $$\Pi(a)\downarrow_{\Pi(\mathcal{B}_\xi\cup\mathcal{C}_\eta)}F.$$
 \end{claim}
\textit{Proof of Claim \ref{claim3final}.}
Since $\alpha$ is $\lambda$-cofinal, $\lambda>2^\omega$, and $|X|=2^\omega$, there is a very good $\beta<\alpha$ such that $X\cap\Gamma_g^{\alpha}\subset \Gamma_g^\beta$. Let $\xi<\eta$ be such that $\mathcal{B}_\xi\subseteq \Gamma_f^\alpha\backslash\Gamma_f^\beta$ and $\xi\notin J_f^\beta$.
It is not difficult to see that 
$\Pi(\mathcal{B}_\xi)\downarrow_{\Pi(\mathcal{A})}D$,
and since $\Pi(\mathcal{C}_\eta)\subseteq D$, $\Pi(\mathcal{B}_\xi)\downarrow_{\Pi(\mathcal{C}_\eta)}D$.
\begin{subclaim}\label{Claim 4.14.3}
There is $a\in I_{\xi\eta}\backslash (I_{\xi\eta}\restriction \omega)$ such that $\Pi(a)\notin E$ and $$\Pi(a)\downarrow_{\Pi(\mathcal{B}_\xi\cup\mathcal{C}_\eta)}E.$$
\end{subclaim}
\textit{Proof of Subclaim \ref{Claim 4.14.3}.}
Suppose, towards a contradiction, that for every $a\in I_{\xi\eta}\backslash (I_{\xi\eta}\restriction\omega)$, $\Pi(a)\not\downarrow_{\Pi(\mathcal{B}_\xi\cup\mathcal{C}_\eta)}E$. Then, for every $a\in I_{\xi\eta}\backslash (I_{\xi\eta}\restriction\omega)$ there is $b_a\in E$ such that $\Pi(a)\not\downarrow_{\Pi(\mathcal{B}_\xi\cup\mathcal{C}_\eta)}b_a$.
The model $E$ was defined as an $a$-primary model over $\Gamma_g^{\alpha+1}\cup X$, therefore $|E|\leq \lambda(T)+(|\Gamma_g^{\alpha+1}\cup X|+\omega)^{<\omega}$. Since $\lambda(T)\leq 2^\omega$ and $|X|=2^\omega$, we obtain $|E|\leq 2^\omega+ |\Gamma_g^{\alpha+1}|$, by Fact \ref{F.4.6}, we get $|E|\leq g(\alpha)$ and $|E|<f(\alpha)$. Since $|I_{\xi\eta}|=f(\alpha)$, then there is $b\in E$ and $J=\{c_i~|~i<\omega\}$, a subset of $I_{\xi\eta}\backslash (I_{\xi\eta}\restriction\omega)$ such that for every $i<\omega$, $\Pi(c_i)\not\downarrow_{\Pi(\mathcal{B}_\xi\cup\mathcal{C}_\eta)}b$ holds. Since $\Pi(I_{\xi\eta}\backslash (I_{\xi\eta}\restriction\omega))$ is independent over $\Pi(\mathcal{B}_\xi\cup\mathcal{C}_\eta)$, then $b\not\downarrow_{\Pi(\mathcal{B}_\xi\cup\mathcal{C}_\eta)\cup\{\Pi(c_j)~|~j<i\}}\Pi(c_i)$ for every $i<\omega$. So $T$ is not superstable, a contradiction. This finishes the proof of Subclaim \ref{Claim 4.14.3}.

Notice that $\Pi(I_{\xi\eta})$ is indiscernible over $\Pi(\mathcal{B}_\xi\cup\mathcal{C}_\eta)$. Since $\Pi(\mathcal{B}_\xi)\downarrow_{\Pi(\mathcal{C}_\eta)}D$, then by domination we get $M_3\downarrow_{\Pi(\mathcal{C}_\eta)}D$, where $M_3$ is an $a$-primary model over $\Pi(\mathcal{B}_\xi\cup\mathcal{C}_\eta)$. So the models $M_0=M_0'=\Pi(\mathcal{A})$, $M_1=M_1'=\Pi(\mathcal{B}_\xi)$, $M_2=\Pi(\mathcal{C}_\eta)$ and $M_2'=D$ satisfy the assumptions of Lemma \ref{L.3.11}, therefore $\Pi(I_{\xi\eta})$ is indiscernible over $\Pi(\mathcal{B}_\xi)\cup D$. By Remark \ref{R.3.12}, if $M'_3$ is an $a$-primary model over $\Pi(\mathcal{B}_\xi)\cup D$ with $\Pi(I_{\xi\eta}\restriction \omega)\subseteq M_3'$, then $Av(\Pi(I_{\xi\eta}\restriction \omega),M_3')\perp D$ and $\Pi(I_{\xi\eta})$ is independent over $\Pi(\mathcal{B}_\xi)\cup D$. So, if $a$ is the element given in Subclaim \ref{Claim 4.14.3} and $\Pi(a)\notin M_3'$ holds, then $tp(\Pi(a),M_3')\perp D$.
\begin{subclaim}\label{Claim 4.14.4}
$tp(\Pi(a),E)\perp D$ 
\end{subclaim}
\textit{Proof of Subclaim \ref{Claim 4.14.4}.}
Let $M_3'$ be an $a$-primary model over $\Pi(\mathcal{B}_\xi)\cup D$ with $\Pi(I_{\xi\eta}\restriction \omega)\subseteq M_3'$. Since $E$ is $a$-saturated, then there is $\mathcal{F}:M_3'\rightarrow E$ an elementary embedding such that $\mathcal{F}\restriction \Pi(\mathcal{B}_\xi)\cup D=id$. Let $b$ be such that $b\models \mathcal{F}(Av(\Pi(I_{\xi\eta}\restriction \omega),M_3'))$, since $Av(\Pi(I_{\xi\eta}\restriction \omega),M_3')\perp D$, then $tp(b,\mathcal{F}(M_3'))\perp D$. By the way $I_{\xi\eta}$ was chosen and Remark \ref{R.3.12}, we know that $\Pi(I_{\xi\eta})$ is independent over $\Pi(\mathcal{B}_\xi)\cup D$, by Lemma \ref{L.3.9} we conclude that $\mathcal{F}(Av(\Pi(I_{\xi\eta}\restriction \omega),M_3'))$ doesn't fork over $\Pi(\mathcal{B}_\xi)\cup D$. On the other hand, by Subclaim \ref{Claim 4.14.3} $\Pi(a)\downarrow_{\Pi(\mathcal{B}_\xi\cup\mathcal{C}_\eta)}E$, so $\Pi(a)\downarrow_{\Pi(\mathcal{B}_\xi)\cup D}\mathcal{F}(M_3')$. By Fact \ref{F.3.5}, since $tp(b,\mathcal{F}(M_3'))\perp D$, $b\downarrow_{\Pi(\mathcal{B}_\xi)\cup D}\mathcal{F}(M_3')$ and $\Pi(a)\downarrow_{\Pi(\mathcal{B}_\xi)\cup D}\mathcal{F}(M_3')$ hold, then $tp(\Pi(a),\mathcal{F}(M_3'))\perp D$.

To show that $tp(\Pi(a),E)\perp D$ let $d$ and $B$ be such that $d\downarrow_DE$, $D\subseteq B$, $\Pi(a)\downarrow_EB$, and $d\downarrow_EB$.
By transitivity, $d\downarrow_DE$ and $d\downarrow_EB$ implies that $d\downarrow_DE\cup B$. By Subclaim \ref{Claim 4.14.3} we know that $\Pi(a)\downarrow_{\Pi(\mathcal{B}_\xi\cup\mathcal{C}_\eta)}E$, then by transitivity we get $\Pi(a)\downarrow_{\Pi(\mathcal{B}_\xi\cup\mathcal{C}_\eta)}E\cup B$. Therefore $d\downarrow_D\mathcal{F}(M_3')\cup B$ and $\Pi(a)\downarrow_{\Pi(\mathcal{B}_\xi)\cup D}\mathcal{F}(M_3')\cup B$ hold, so $d\downarrow_D\mathcal{F}(M_3')$, $d\downarrow_{\mathcal{F}(M_3')} B$ and $\Pi(a)\downarrow_{\mathcal{F}(M_3')} B$ hold. Since $tp(\Pi(a),\mathcal{F}(M_3'))\perp D$, we conclude that $\Pi(a)\downarrow_Bb$, finishing the proof of Subclaim \ref{Claim 4.14.4}.

Let $I_X$ be the set $\bigcup\{\mathcal{B}_r, I_{rp}~|~\mathcal{B}_r\not\subseteq\Gamma_g^{\alpha+1}\wedge r<p\wedge\mathcal{C}_p\subseteq X\backslash \Gamma_g^{\alpha+1}\}$. Let us show that $D\downarrow_XI_X\cup \Gamma_g^{\alpha+1}$.

If $D\not\downarrow_{X}I_X\cup\Gamma_g^{\alpha+1}$, then there are finite $c\in D$ and $b\in (I_X\cup\Gamma_g^{\alpha})\backslash X$ such that $a\not\downarrow_{X}b$.
Since $D$ is $a$-constructable over $X$, then it is $a$-atomic over $X$. So, there is a finite $A_1\subseteq X$ such that $stp(c,A_1)\vdash tp(c,X)$.
Since $T$ is superstable, there is a finite $A_2\subseteq X$ such that $c\cup b\downarrow _{A_2}X$. Denote by $A$ the set $A_1\cup A_2$. Since $X$ is nice, $A$ is a finite subset of $X$ and $b\in (I_X\cup\Gamma_g^{\alpha})\backslash X$, then by Lemma \ref{L.4.10} there is $\mathcal{F}\in Saut(\mathcal{M},A)$ such that $\mathcal{F}(b)\in X$. Therefore $stp(\mathcal{F}(c),A_1)\vdash tp(c,X)$, and $\mathcal{F}(c)\downarrow_{A_1}X$, we conclude $\mathcal{F}(c)\downarrow_{A}\mathcal{F}(b)$ and $c\downarrow_{A}b$. Since $c\cup b\downarrow _{A_2}X$, then $c\cup b\downarrow _{A}X$. Therefore $c\downarrow_{X}b$, a contradiction.

By Fact \ref{F.4.12}, we know that $I_X\cup X\downarrow_{X\cap \Gamma_g^{\alpha+1}}\Gamma_g^{\alpha+1}$, then $I_X\downarrow_{X}\Gamma_g^{\alpha+1}$. Since $D\downarrow_XI_X\cup \Gamma_g^{\alpha+1}$, we conclude that $I_X\downarrow_D\Gamma_g^{\alpha+1}$.
By the way $E$ was chosen, we know that $E$ is $a$-constructible over $D\cup \Gamma_g^{\alpha+1}$. Since $D$ is $a$-saturated, we get that $\Gamma_g^{\alpha+1}\rhd_DE$. By domination we conclude $I_X\downarrow_DE$.
Therefore, for every $c\in I_X$ we have that $c\downarrow_DE$. Since $c\downarrow_EE$ and $\Pi(a)\downarrow_EE$ hold, then by Subclaim \ref{Claim 4.14.4} we conclude that $c\downarrow_E\Pi(a)$ for every $c\in I_X$. By the finite character we get $I_X\downarrow_E\Pi(a).$
By the way $F$ was chosen, we know that $F$ is $a$-constructible over $I_X\cup E$, and since $E$ is $a$-saturated, we conclude that $I_X\rhd_EF.$ Therefore $F\downarrow_E\Pi(a)$.
Since $\Pi(a)\downarrow_{\Pi(\mathcal{B}_\xi\cup\mathcal{C}_\eta)}E$, by transitivity we conclude $\Pi(a)\downarrow_{\Pi(\mathcal{B}_\xi\cup\mathcal{C}_\eta)}F$. This finishes the proof of Claim \ref{claim3final}.

\begin{claim}\label{claim4final}
There is an independent indiscernible sequence $(b_i)_{i<f(\alpha)^+}$ over $\Pi(\mathcal{B}_\xi\cup\mathcal{C}_\eta)$ in $\mathcal{A}^g$.
\end{claim}
\textit{Proof of Claim \ref{claim4final}}
Recall that $\Pi(a)\in \mathcal{A}^g$ and $\mathcal{A}^g$ is $a$-constructable over $F\cup \Gamma_g$. Then $\mathcal{A}^g$ is $a$-atomic over $F\cup \Gamma_g$ and there is a finite $B\subseteq F\cup \Gamma_g$ such that $(tp(\Pi(a),F\cup \Gamma_g),B)\in F_\omega^a$ (it is $a$-isolated) and $\Pi(a)\in \mathcal{N}$, where $\mathcal{N}\subseteq \mathcal{A}^g$ is $a$-primary over $F\cup B$. Let $B'=B\backslash F$, there is a nice set $\mathcal{Y}$ such that $\mathcal{Y}\cap F=\mathcal{A}$, $B'\subseteq \mathcal{Y}$, $\mathcal{Y}$ $\Gamma_g$-nice for all $(r,p)$ that satisfy $\mathcal{B}_r,\mathcal{C}_p\subset \mathcal{Y}$, and $S=\{r\in J_g\mid (r\in (J_g)_{<\lambda}\wedge \mathcal{B}_r\subset \mathcal{Y})\vee(r\in (J_g)_{\lambda}\wedge \mathcal{C}_r\subset \mathcal{Y})\}$ is finite. Define $\mathcal{X}=\{r\in J_g\mid (r\in (J_g)_{<\lambda}\wedge \mathcal{B}_r\subset X)\vee(r\in (J_g)_{\lambda}\wedge \mathcal{C}_r\subset X)\}$. Let $\bar{S}=S\cup \{r\in (J_g)_{<\lambda}\mid\exists p\in S\ (r<p)\}$ and $\bar{\mathcal{X}}=\mathcal{X}\cup \{r\in (J_g)_{<\lambda}\mid\exists p\in \mathcal{X}\ (r<p)\}$. By the way $\bar{\mathcal{X}}$ was defined, we know that for every limit ordinal $\theta<\lambda$ and $\zeta\in J_g$, if for all $\theta'<\theta$, $\zeta\restriction{\theta'}\in \bar{\mathcal{X}}$ holds, then $\zeta\restriction{\theta}\in \bar{\mathcal{X}}$. Notice that since $cf(\alpha)=\lambda$, if $\theta<\lambda$ is a limit ordinal such that for all $\theta'<\theta$, $\zeta\restriction{\theta'}\in J_g^{\alpha+1}$ holds, then $\zeta\restriction{\theta}\in J_g^{\alpha+1}$. We conclude that if $\theta<\lambda$ and $\zeta\in J_g$ are such that for all $\theta'<\theta$, $\zeta\restriction\theta'\in \bar{\mathcal{X}}\cup J_g^{\alpha+1}$ and $\zeta\restriction{\theta}\in \bar{S}\backslash (\bar{\mathcal{X}}\cup J_g^{\alpha+1})$, then $\theta$ is a successor ordinal.
Let $\{u_i\}_{i<f(\alpha)^+}$ be a sequence of subtrees of $J_g$ with the following properties:
\begin{itemize}
\item $u_0=\bar{S}$.
\item Every $u_i$ is a tree isomorphic to $u_0$.
\item If $i\neq j$, then $u_i\cap u_j=u_0\cap (\bar{\mathcal{X}}\cup J_g^{\alpha+1})$.
\item Every $\zeta\in dom(c_g)\cap u_0$ satisfies $c_f(\zeta)=c_f(G_i( \zeta))$, where $G_i$ is the isomorphism between $u_0$ and $u_i$.
\end{itemize}
For every $\zeta\in u_0$ and $\theta<\lambda$ such that $\zeta\restriction\theta\in \bar{\mathcal{X}}\cup J_g^{\alpha+1}$ and $\zeta\restriction{\theta+1}\in u_0\backslash (\bar{\mathcal{X}}\cup J_g^{\alpha+1})$, it holds by Definition \ref{D.2.6} that $\zeta\restriction\theta$ has $\kappa$ many immediate successors in $J_g\backslash J_g^{\alpha+1}$. Also by Definition \ref{D.2.6} the elements of $J_f$ are all the functions $\eta:s\rightarrow \lambda\times \kappa^4$ that satisfy the items 1 to 8, therefore each of the immediate successors of $\zeta\restriction\gamma$, $\zeta'$, satisfies that in the set $\{r\in J_f~|~\zeta'\leq r\}$ there is a subtree isomorphic (as coloured tree) to $\{p\in u_0\backslash (\bar{\mathcal{X}}\cup J_g^{\alpha+1})\mid \zeta\restriction \gamma+1\leq p\}$.
This and the fact that $S$ is finite, gives the existence of the sequence $\{u_i\}_{i<f(\alpha)^+}$.
By the way we chose the sequence $\{u_i\}_{i<f(\alpha)^+}$, for every $i<f(\alpha)^+$, the isomorphism $G_i$ induces a coloured trees isomorphism $\bar{G}_i:\bar{\mathcal{X}}\cup J_g^{\alpha+1}\cup u_0\rightarrow \bar{\mathcal{X}}\cup J_g^{\alpha+1}\cup u_i$ such that $\bar{G}_i\restriction \bar{\mathcal{X}}\cup J_g^{\alpha+1}= id$. Let us denote by $z_i$ the tree $\bar{\mathcal{X}}\cup J_g^{\alpha+1}\cup u_i$.

Let us define $U_i=\{\mathcal{B}_r\mid r\in z_i\wedge r\in (J_g)_{<\lambda}\}\cup\{\mathcal{C}_p\mid p\in z_i\wedge p\in (J_g)_{\lambda}\}$ and $\bar{U}_i=U_i\cup \{I_{rp}\mid \mathcal{B}_r\in U_i\wedge \mathcal{C}_p\in U_i\wedge r<p\}$. Notice that $\bigcup \bar{U}_i$ is nice for all $i<f(\alpha)^+$. Since $u_i$ is isomorphic to $\bar{S}$, then $p\in z_i$ and $r<p$, implies $r\in z_i$. Therefore, $\bigcup\bigcup_{j\neq i}\bar{U}_j$ is nice for all $i<f(\alpha)^+$.
\begin{subclaim}\label{Claim 4.14.5}
For all $i<f(\alpha)^+$ it holds that $\bigcup \bar{U}_i\downarrow_F\bigcup\bigcup_{j\neq i}\bar{U}_j$.
\end{subclaim}
\textit{Proof of Subclaim \ref{Claim 4.14.5}.}
By the way the sets $\bar{U}_i$ were constructed, we know that $(\bigcup \bar{U}_i) \cap(\bigcup \bar{U}_j)=\Gamma_g^{\alpha+1}\cup X\cup I_X$ for all $i\neq j$. Let us denote by $\mathbb{F}$ the set $\Gamma_g^{\alpha+1}\cup X\cup I_X$. By Corollary 4.13 we know that $$\bigcup \bar{U}_i\downarrow_\mathbb{F}\bigcup\bigcup_{j\neq i}\bar{U}_j.$$
Let us prove that $F\downarrow_\mathbb{F}\bigcup\bigcup_{j<f(\alpha)^+}\bar{U}_j$. Suppose it is false, then $F\not\downarrow_\mathbb{F}\bigcup\bigcup_{j<f(\alpha)^+}\bar{U}_j$ and there are finite $c\in F$ and $b\in \bigcup\bigcup_{j<f(\alpha)^+}\bar{U}_j$ such that $c\not\downarrow_{\mathbb{F}}b$.
Since $F$ is $a$-constructable over $\mathbb{F}$, then it is $a$-atomic over $\mathbb{F}$. So, there is a finite $A_1\subseteq \mathbb{F}$ such that $stp(c,A_1)\vdash tp(c,\mathbb{F})$.
Since $T$ is superstable, there is a finite $A_2\subseteq \mathbb{F}$ such that $c\cup b\downarrow _{A_2}\mathbb{F}$. Denote by $A$ the set $A_1\cup A_2$. By Lemma \ref{L.4.10} there is $\mathcal{F}\in Saut(\mathcal{M},A)$ such that $\mathcal{F}(b)\in \mathbb{F}$. Therefore $stp(\mathcal{F}(c),A_1)\vdash tp(c,\mathbb{F})$, and $\mathcal{F}(c)\downarrow_{A_1}\mathbb{F}$. So $\mathcal{F}(c)\downarrow_{A}\mathcal{F}(b)$ and $c\downarrow_{A}b$. Since $c\cup b\downarrow _{A_2}\mathbb{F}$, then $c\cup b\downarrow _{A}\mathbb{F}$. Therefore $c\downarrow_{\mathbb{F}}b$, a contradiction.

Since $F\downarrow_\mathbb{F}\bigcup\bigcup_{j<f(\alpha)^+}\bar{U}_j$ and $\bigcup \bar{U}_i\downarrow_\mathbb{F}\bigcup\bigcup_{j\neq i}\bar{U}_j$ holds, we conclude that $\bigcup \bar{U}_i\downarrow_F\bigcup\bigcup_{j\neq i}\bar{U}_j$, finishing the proof of Subclaim \ref{Claim 4.14.5}.

The isomorphisms $(\bar{G}_i)_{i<f(\alpha)^+}$ induce the following elementary maps $\mathcal{H}^i_{rp}:\mathcal{B}_r\cup\mathcal{C}_p\cup I_{rp}\rightarrow\mathcal{B}_{\bar{G}_i(r)}\cup\mathcal{C}_{\bar{G}_i(p)}\cup I_{\bar{G}_i(r)\bar{G}_i(p)}$ for all $r,p\in z_0$ ($r\in (J_g)_{<\lambda}$ and $p\in (J_g)_{\lambda}$), given by $\mathcal{H}^i_{rp}=H_{\bar{G}_i(r)\bar{G}_i(p)}\circ H^{-1}_{rp}$. 

It is easy to check that the map $\mathcal{H}_i:\bigcup\bar{U}_0\rightarrow\bigcup\bar{U}_i$ defined by $$\mathcal{H}_i=\bigcup_{\eta\in z_0\cap (J_f)_\lambda} ~ \bigcup_{\xi\in z_0\cap (J_f)_{<\lambda}, \xi<\eta}\mathcal{H}^i_{\xi\eta}$$ is elementary.
Notice that for any permutation $\mathcal{R}: f(\alpha)^+\rightarrow f(\alpha)^+$ and any $i<f(\alpha)^+$, $tp(\bigcup\bigcup_{j<i}\bar{U}_j,\Gamma_g^{\alpha+1}\cup X\cup I_X)=tp(\bigcup\bigcup_{j<i}\bar{U}_{\mathcal{R}(j)},\Gamma_g^{\alpha+1}\cup X\cup I_X)$ holds.

Therefore $(\bigcup\bar{U}_i)_{i<f(\alpha)^+}$ is an indiscernible sequence over $\Gamma_g^{\alpha+1}\cup X\cup I_X$. So, for all $i<f(\alpha)^+$, $stp(\bigcup\bar{U}_0,\Gamma_g^{\alpha+1}\cup X\cup I_X)=stp(\bigcup\bar{U}_i,\Gamma_g^{\alpha+1}\cup X\cup I_X)$. Let $\mathcal{G}_i: F\cup\bigcup\bar{U}_0\rightarrow F\cup\bigcup\bar{U}_i$, be given by $\mathcal{G}_i\restriction F=id$ and $\mathcal{G}_i\restriction \bigcup\bar{U}_0=\mathcal{H}_i$. 
It is easy to check that $\mathcal{G}_i$ is elementary.

Let us define for all $i<f(\alpha)^+$ the model $M_i\subseteq \mathcal{A}^g$ as an $a$-primary model over $F\cup\bigcup_{j<i}M_j\cup\bigcup\bar{U}_i$, with $\mathcal{N}\subseteq M_0$ and let $b_0\in M_0$ be $\Pi(a)$ (notice that $B\subseteq \bar{U}_0$ was chosen such that $(tp(\Pi(a),F\cup \Gamma_g),B)\in F_\omega^a$ and $\Pi(a)\in \mathcal{N}$, $\mathcal{N}$ is the $a$-primary model over $F\cup B$). For all $0<i<f(\alpha)^+$ let $\bar{\mathcal{G}}_i\in Saut(\mathcal{M},\Gamma_g^{\alpha+1}\cup X\cup I_X)$ be such that $\bar{\mathcal{G}}_i\restriction F\cup\bigcup\bar{U}_i=\mathcal{G}_i\restriction F\cup\bigcup\bar{U}_i$ and $b_i\in M_i$ be such that $stp(b_i,\mathcal{G}_i(B))=stp(\bar{\mathcal{G}}_i(\Pi(a)),\mathcal{G}_i(B))$. We know that $(tp(\Pi(a),F\cup \Gamma_g),B)\in F_\omega^a$, so by $a$-isolation and the definition of $\bar{\mathcal{G}}_i$ we conclude that $(tp(b_i,\bar{\mathcal{G}}_i(F \cup\bigcup\bar{U}_0)),\mathcal{G}_i(B))\in F_\omega^a$, so $(tp(b_i,F \cup\bigcup\bar{U}_i),\mathcal{G}_i(B))\in F_\omega^a$. Therefore $tp(b_i,F)=tp(\bar{\mathcal{G}}_i(\Pi(a)),F)$ and since $\bar{\mathcal{G}}_i$ is an automorphism that fix $F$, we conclude that $tp(b_i,F)=tp(\Pi(a),F)$. On the other hand $(tp(b_i,F \cup\bigcup\bar{U}_i),\mathcal{G}_i(B))\in F_\omega^a$ implies that $b_i\cup F \cup\bigcup\bar{U}_i$ is $a$-constructable over $F \cup\bigcup\bar{U}_i$, since $F$ is $a$-saturated then $\bigcup\bar{U}_i\rhd_Fb_i\cup\bigcup\bar{U}_i$. By Subclaim \ref{Claim 4.14.5} we know that $\bigcup \bar{U}_i\downarrow_F\bigcup\bigcup_{j\neq i}\bar{U}_j$, so by domination we conclude that $b_i\cup\bigcup \bar{U}_i\downarrow_F\bigcup\bigcup_{j\neq i}\bar{U}_j$, in particular $b_i\downarrow_F\bigcup\bigcup_{j\neq i}\bar{U}_j$ holds for all $i<f(\alpha)^+$. 

Notice that for all $i<f(\alpha)^+$, $M_i$ is $a$-constructable over $F\cup\bigcup\bigcup_{j\leq i}\bar{U}_j$.
Therefore $\bigcup\bigcup_{k\leq j}\bar{U}_k\rhd_FM_j$ holds for all $i<f(\alpha)^+$, and since $b_i\downarrow_F\bigcup\bigcup_{j\neq i}\bar{U}_j$ holds for all $i<f(\alpha)^+$, then $b_i\downarrow_FM_j$ holds for all $j,i<f(\alpha)^+$, $j<i$. In particular $b_i\downarrow_F\bigcup_{k\leq j}b_k$ holds for all $j,i<f(\alpha)^+$, $j<i$. We conclude that $b_i\downarrow_F\bigcup_{j< i}b_j$ holds for all $i<f(\alpha)^+$.
Since $tp(b_i,F)=tp(\Pi(a),F)$ and $\Pi(a)\downarrow_{\Pi(\mathcal{B}_\xi\cup\mathcal{C}_\eta)}F$, we get that $b_i\downarrow_{\Pi(\mathcal{B}_\xi\cup\mathcal{C}_\eta)}F$ and by transitivity we conclude that $b_i\downarrow_{\Pi(\mathcal{B}_\xi\cup\mathcal{C}_\eta)}\bigcup_{j< i}b_j$. So $(b_i)_{i<f(\alpha)^+}$ is an independent sequence over $\Pi(\mathcal{B}_\xi\cup\mathcal{C}_\eta)$. Since for $i\neq j$ we know that $tp(b_i,F)=tp(b_j,F)$, the types over $F$ are stationary, and $b_i\downarrow_F\bigcup_{j< i}b_j$, then we conclude that $(b_i)_{i<f(\alpha)^+}$ is an indiscernible sequence over $F$. This finishes the proof of Claim \ref{claim4final}.

For every $i<f(\alpha)^+$ let $c_i$ be $\Pi^{-1}(b_i)$, since $\Pi$ is an isomorphism, then $(c_i)_{i<f(\alpha)^+}$ is an indiscernible sequence over $\mathcal{B}_\xi\cup\mathcal{C}_\eta$ and an independent sequence over $\mathcal{B}_\xi\cup\mathcal{C}_\eta$, notice that $c_0=a$, so $c_0\in I_{\xi\eta}$. This contradicts Lemma \ref{lemma_no_sequence}.

\end{proof}

\section{Corollaries}
Recall that $\lambda$ is the cardinal $(2^\omega)^+$.
\begin{corollary}\label{C.4.15}
If $\kappa$ is inaccessible, and $T$ is a theory with S-DOP, then $E_{\lambda\text{-club}}^\kappa\leq_c\cong_T$.
\end{corollary}
\begin{proof}
Let $f$ and $g$ be elements of $\kappa^\kappa$.
First we will construct a function $F:\kappa^\kappa\rightarrow\kappa^\kappa$ such that $f\ E^\kappa_{\lambda\text{-club}}\ g$ if and only if $\mathcal{A}^{F(f)}$ and $\mathcal{A}^{F(g)}$ are isomorphic. This is to get over the restrictions on $f$ and $g$ in the hypothesis of Theorem \ref{T.4.14} .

For every cardinal $\alpha<\kappa$, define $S_\alpha=\{\beta\in Card\cap\kappa~|~\lambda,\alpha^{+++},\alpha^\lambda<\beta\}$.
Let $\mathcal{G}_\beta$ be a bijection from $\kappa$ into $S_\beta$, for every $\beta<\kappa$. For every $f\in \kappa^\kappa$ define $F(f)$ by $F(f)(\beta)=\mathcal{G}_\beta(f(\beta))$, for every $\beta<\kappa$. Clearly $f\ E^\kappa_{\lambda\text{-club}}\ g$ if and only if $F(f)\  E^\kappa_{\lambda\text{-club}}\ F(g)$ i.e. $\mathcal{A}^{F(f)}$ and $\mathcal{A}^{F(g)}$ are isomorphic and $F$ is continuous.

Finally we need to find $\mathcal{G}:\{F(f)~|~f\in \kappa^\kappa\}\rightarrow\kappa^\kappa$ such that $\mathcal{A}_{\mathcal{G}(F(f))}\cong \mathcal{A}^{F(f)}$ and $f \mapsto \mathcal{G}(F(f))$ is continuous.
Notice that for every $f,g\in \kappa^\kappa$ and $\alpha<\kappa$, by Definition \ref{D.2.6} and the definition of $J_f^\alpha$ in Remark \ref{R.2.7a}, it holds: $$F(f)\restriction\alpha=F(g)\restriction\alpha\Leftrightarrow J_{F(f)}^\alpha=J_{F(g)}^\alpha.$$
By Definition \ref{D.4.5}, for every $f,g\in \kappa^\kappa$ and $\alpha<\kappa$ it holds:
$$J_{F(f)}^\alpha=J_{F(g)}^\alpha\Leftrightarrow \Gamma_{F(f)}^\alpha=\Gamma_{F(g)}^\alpha.$$
By the definition of $\mathcal{A}_f^\alpha$ in Theorem \ref{T.4.14}, for every $f,g\in \kappa^\kappa$ and $\alpha<\kappa$ an $F(f)$-good and $F(g)$-good cardinal, it holds:
$$\Gamma_{F(f)}^\alpha=\Gamma_{F(g)}^\alpha\Leftrightarrow \mathcal{A}_{F(f)}^\alpha\cong\mathcal{A}_{F(g)}^\alpha.$$
In general,since there are club many $F(f)$-good and $F(g)$-good cardinals, then by the definition of $\mathcal{A}_f^\alpha$ in Theorem \ref{T.4.14} we can construct the models $\mathcal{A}^f$ such that for every $f,g\in \kappa^\kappa$ and $\alpha<\kappa$, it holds: $$J_{F(f)}^\alpha=J_{F(g)}^\alpha\Leftrightarrow \mathcal{A}_{F(f)}^\alpha=\mathcal{A}_{F(g)}^\alpha.$$
So we can construct the models $\mathcal{A}^f$ such that for every $f,g\in \kappa^\kappa$ and $\alpha<\kappa$, it holds: $$F(f)\restriction\alpha=F(g)\restriction\alpha\Leftrightarrow \mathcal{A}_{F(f)}^\alpha=\mathcal{A}_{F(g)}^\alpha.$$
For every $f\in \kappa^\kappa$ define $C_f\subseteq Card\cap\kappa$ such that $\forall \alpha\in C_f$, it holds that for all $\beta$ ordinal smaller than $\alpha$, $\mid\mathcal{A}^\beta_{F(f)}\mid<\mid\mathcal{A}^\alpha_{F(f)}\mid$. For every $f\in \kappa^\kappa$ and $\alpha\in C_f$ choose $E^\alpha_f:dom(\mathcal{A}^\alpha_{F(f)})\rightarrow\mid\mathcal{A}^\alpha_{F(f)}\mid$ a bijection, such that $\forall \beta,\alpha\in C_f$, $\beta<\alpha$ it holds that $E_f^\beta\subseteq E_f^\alpha$. Therefore $\bigcup_{\alpha\in C_f}E_f^\alpha=E_f$ is such that $E_f:dom(\mathcal{A}^{F(f)})\rightarrow\kappa$ is a bijection, and for every $f,g\in \kappa^\kappa$ and $\alpha<\kappa$ it holds:
If $F(f)\restriction\alpha=F(g)\restriction\alpha$, then $E_f\restriction dom(\mathcal{A}_{F(f)}^\alpha)=E_g\restriction dom(\mathcal{A}_{F(g)}^\alpha)$.

Let $\pi$ be the bijection in Definition \ref{D.1.6}, define the function $\mathcal{G}$ by: $$\mathcal{G}(F(f))(\alpha)=\begin{cases} 1 &\mbox{if } \alpha=\pi(m,a_1,a_2,\ldots,a_n) \text{ and } \mathcal{A}^{F(f)}\models P_m(E_f^{-1}(a_1),E_f^{-1}(a_2),\ldots,E_f^{-1}(a_n))\\
0 & \mbox{in other case. } \end{cases}$$
To show that $\mathcal{G}$ is continuous, let $[\eta\restriction\alpha]$ be a basic open set and $\xi\in \mathcal{G}^{-1}[[\eta\restriction\alpha]]$. So, there is $\beta\in C_\xi$ such that for all $\gamma<\alpha$, if $\gamma=\pi(m,a_1,a_2,\ldots,a_n)$, then $E^{-1}_\xi(a_i)\in dom(\mathcal{A}_\xi^\beta)$ holds for all $i\leq n$. Since for all $\zeta\in [\xi\restriction\beta]$ it holds that $\mathcal{A}_\xi^\beta=\mathcal{A}_\zeta^\beta$, then for every $\gamma<\alpha$ that satisfies $\gamma=\pi(m,a_1,a_2,\ldots,a_n)$, it holds that: $$\mathcal{A}^{\xi}\models P_m(E_\xi^{-1}(a_1),E_\xi^{-1}(a_2),\ldots,E_\xi^{-1}(a_n))\Leftrightarrow\mathcal{A}^{\zeta}\models P_m(E_\zeta^{-1}(a_1),E_\zeta^{-1}(a_2),\ldots,E_\zeta^{-1}(a_n)).$$
We conclude that $\mathcal{G}(\zeta)\in [\eta\restriction\alpha] $, and $\mathcal{G}$ is continuous.
\end{proof}

In \cite{HM} it was proved that if $T$ is a classifiable theory and $\mu<\kappa$ is a regular cardinal, then $\cong_T$ is continuously reducible to $E_{\mu\textit{-club}}^\kappa$.

\begin{corollary}\label{C.4.16}
If $\kappa$ is an inaccessible and $T_1$ is a classifiable theory and $T_2$ is a superstable theory with S-DOP, then $\cong_{T_1}\leq_c\cong_{T_2}$.
\end{corollary}

The last corollaries are about $\Sigma_1^1$-completeness.
Suppose $E$ is an equivalence relation on $\kappa^\kappa$. We say that $E$ is $\Sigma_1^1$
  if $E$ is the projection of a closed set in $\kappa^\kappa\times \kappa^\kappa\times \kappa^\kappa$ and it 
  is $\Sigma_1^1$-complete if every $\Sigma_1^1$ equivalence relation is 
  Borel reducible to $E$.

In \cite{HK} it was proved,  under the assumption $V=L$, that $E^\kappa_{\mu\text{-club}}$ is $\Sigma_1^ 1$-complete for all regular $\mu<\kappa$. In \cite{FMR}, under the assumption GCH, it was proved that there exists a cofinality-preserving GCH-preserving forcing extension in which $E^\kappa_{\mu\text{-club}}$ is $\Sigma_1^ 1$-complete for all regular $\mu<\kappa$.

\begin{corollary}\label{C.4.19}
\begin{itemize}
\item Suppose $V=L$. If $\kappa$ is an inaccessible and $T$ is a superstable theory with S-DOP, then $\cong_{T}$ is $\Sigma_1^ 1$-complete.
\item Suppose GCH. There exists a cofinality-preserving GCH-preserving forcing extension in which If $\kappa$ is an inaccessible and $T$ is a superstable theory with S-DOP, then $\cong_{T}$ is $\Sigma_1^ 1$-complete.
\end{itemize}

\end{corollary}

\providecommand{\bysame}{\leavevmode\hbox to3em{\hrulefill}\thinspace}
\providecommand{\MR}{\relax\ifhmode\unskip\space\fi MR }
\providecommand{\MRhref}[2]{%
  \href{http://www.ams.org/mathscinet-getitem?mr=#1}{#2}
}
\providecommand{\href}[2]{#2}

\end{document}